\documentclass[12pt,a4paper,fleqn]{article}
\usepackage{a4wide,amsfonts,amsmath,latexsym,amssymb,euscript,graphicx,units,mathrsfs}

\usepackage{graphicx}
\usepackage{color}
\usepackage{amssymb}
\usepackage[T1]{fontenc}
\usepackage{latexsym}
\usepackage{xypic}
\usepackage{eufrak}
\usepackage{euscript}
\usepackage{amsfonts,amsmath}
\usepackage{verbatim}
\usepackage{fancyhdr}
\usepackage[english]{babel}
\usepackage{mathrsfs}
\usepackage{units}

\newtheorem{prop}{Proposition}[section]

\newtheorem{lem}[prop]{Lemma}
\newtheorem{rem}[prop]{Remark}

\newtheorem{thm}[prop]{Theorem}
\newtheorem{theo}[prop]{Theorem}

\renewcommand{\geq}{\geqslant}
\def\leq{\leqslant}
\newcommand{\D}{\mathbb{D}}
\newcommand{\N}{\mathbb{N}}

\newcommand{\R}{\mathbb{R}}

\newcommand{\bP}{P}

\def\esp{E}
\def\var{{\mathbb{Var}}}
\def\HH{\EuFrak H}

\def\e{\varepsilon}

\def\1{{\mathbf{1}}}

\def\1{{\mathbf{1}}}
\def\0.5{{\frac{1}{2}}}
\def\var{{\rm{Var}}}

\newcommand{\qed}{\nopagebreak\hspace*{\fill}
{\vrule width6pt height6ptdepth0pt}\par}

\newcounter{rea}
\begin{document}

\begin{center}
{\Large{\bf Convergence in total variation on Wiener chaos}}
\normalsize
\\~\\ by
Ivan Nourdin\footnote{Email: {\tt inourdin@gmail.com}; IN was partially supported by the
ANR grants ANR-09-BLAN-0114 and ANR-10-BLAN-0121.} and
Guillaume Poly\footnote{Email: {\tt guillaume.poly@crans.org}\\}\\
{\it Universit\'e de Lorraine and Universit\'e Paris Est}\\~\\
\end{center}

{\small \noindent {{\bf Abstract}:
Let $\{F_n\}$ be a sequence of random variables belonging to a finite sum of Wiener chaoses.
Assume further that it converges in distribution towards $F_\infty$ satisfying ${\rm Var}(F_\infty)>0$.
Our first result is a sequential version of a theorem by Shigekawa \cite{Shigekawa}. More precisely,
 we prove, without additional assumptions, that the sequence
 $\{F_n\}$ actually converges in total variation and that the law of
 $F_\infty$ is absolutely continuous.
 We give an application to discrete non-Gaussian chaoses.
 In a second part, we assume that each $F_n$ has more specifically the form of a multiple Wiener-It\^o integral
 (of a fixed order) and that it converges in $L^2(\Omega)$ towards $F_\infty$. We then give an upper bound
 for the distance in total variation between the laws of $F_n$ and $F_\infty$.
 As such, we recover an inequality due to Davydov and Martynova \cite{DM}; our rate is weaker compared to
 \cite{DM} (by a power of 1/2), but the advantage is that
 our proof is not only sketched as in \cite{DM}.
Finally, in a third part we show that the convergence in the celebrated Peccati-Tudor theorem actually
 holds in the total variation topology.  


\section{Introduction}
In a seminal paper of 2005, Nualart and Peccati \cite{nunugio} discovered the surprising fact that
convergence in distribution for sequences of multiple Wiener-It\^o integrals to the Gaussian is equivalent to convergence of just the fourth moment.
A new line of research was born.
Indeed, since the publication of this important paper, many improvements and developments on this theme have been considered. (For an overview of the existing literature, we refer the reader to the book \cite{NouPecBook}, to the survey \cite{lecturenotes-coursfondation} or to the constantly updated web page \cite{webpage}.)

Let us only state one of these results, whose proof relies on the combination of Malliavin calculus and Stein's method (see, e.g., \cite[Theorem 5.2.6]{NouPecBook}). When $F,G$ are random variables, we write $d_{TV}(F,G)$ to
indicate the total variation distance between the laws of $F$ and $G$, that is,
\[
d_{TV}(F,G)=\sup_{A\in\mathcal{B}(\R)}\left|\bP(F\in A)-\bP(G\in A)\right|=\frac12\sup_\phi \big|E[\phi(F)]-E[\phi(G)]\big|,
\]
where the first (resp. second) supremum is taken\footnote{One can actually restrict to {\it bounded} Borel sets
without changing the value of the supremum; this easy remark is going to be used many times in the forthcoming proofs.} over Borel sets $A$ of $\R$ (resp. over continuous functions $\phi:\R\to\R$ which are bounded by 1).
\begin{theo}\label{N-P}
If $k\geq 2$ is an integer, if $F$ is an element of the $k$th Wiener chaos $\mathcal{H}_k$ satisfying $E[F^2]=1$ and if $N\sim\mathcal{N}(0,1)$, then
\[
d_{TV}(F,N)\leq \sqrt{\frac{4k-4}{3k}}\sqrt{\left|\esp\left[F^4\right]-3\right|}.
\]
\end{theo}

\bigskip

As an almost immediate corollary of Theorem \ref{N-P}, we get the surprising fact that if a sequence of multiple Wiener-It\^o integrals with unit variance converges in distribution to the standard Gaussian law, then it automatically converges in total variation (\cite[Corollary 5.2.8]{NouPecBook}).
The main thread of the present paper is the seek for other instances where such a phenomenon could occur.
In particular, a pivotal role will be played by the sequences having the form of a (vector of) multiple Wiener-It\^o integral(s) or, more generally, belonging to a {\it finite} sum of Wiener chaoses.
As we said, the proof of Theorem \ref{N-P} relies in a crucial way to the use of Stein's method. In a non-discrete
framework (which is the case here), it is fairly understood that
this method can give good results with respect to the total variation distance only in dimension one (see \cite{ChatterjeeMeckes}) and when
the target law is Gaussian (see \cite{ChenGoldsteinShao}). Therefore, to reach our goal
we need to introduce completely new ideas with respect to the existing literature.
As anticipated, we will manage to exhibit three different situations where the convergence in distribution turns out to be equivalent to the convergence in total variation.
In our new approach, an important role is played by the fact that the Wiener chaoses enjoy many nice properties, such as hypercontractivity (Theorem \ref{hyper}), product formula (\ref{multiplication}) or Hermite polynomial representation of multiple integrals (\ref{hq}).

\bigskip
Let us now describe our main results in more detail.
Our first example focuses on sequences belonging to a
finite sum of chaoses and may be seen as a sequential version of a theorem by Shigekawa \cite{Shigekawa}. More specifically, let $\{F_n\}$ be a sequence in $\bigoplus_{k=0}^p\mathcal{H}_k$
(where $\mathcal{H}_k$ stands for the $k$th Wiener chaos;
by convention $\mathcal{H}_0=\R$), and assume that it
converges in distribution towards a random variable $F_{\infty}$.
Assume moreover that the variance of $F_\infty$ is not zero.
Let $d_{FM}$ denote the Fortet-Mourier distance, defined by
\[
d_{FM}(F,G)=\sup_{\phi} \big|E[\phi(F)]-E[\phi(G)]\big|,
\]
where the supremum is taken over 1-Lipschitz functions $\phi:\R\to\R$ which are bounded by 1.
We prove that there exists a constant $c>0$ such that, for any $n\geq 1$,
\begin{equation}\label{ineq}
d_{TV}(F_n,F_\infty)\leq c\, d_{FM}(F_n,F_\infty)^{\frac{1}{2p+1}}.
\end{equation}
Since it is well-known that $d_{FM}$ metrizes the convergence in distribution (see, e.g., \cite[Theorem 11.3.3]{Dudley book}), our inequality
(\ref{ineq}) implies in particular that $F_n$ converges to $F_\infty$ not only in distribution, but also in total variation.
Besides, one can further prove that the law of $F_\infty$ is absolutely continuous with respect to the Lebesgue measure. This fact is an interesting first step towards a full description of the closure in distribution\footnote{It is worthwhile noting that the Wiener chaoses are closed for the convergence in {\it probability}, as shown by
Schreiber \cite{Sch} in 1969.}
of the Wiener chaoses $\mathcal{H}_k$, which is still an open problem except when $k=1$ (trivial) or $k=2$ (see \cite{G.I}).
We believe that our method is robust enough to be applied to some more general situations, and
here is a short list of possible extensions of (\ref{ineq}) that we plan to study in some
subsequent papers: 
\begin{enumerate}
\item[(i)] extension to the multidimensional case;
\item[(ii)] improvement of the rate of convergence; 
\item[(iii)] extension to other types of chaoses (in the spirit of \cite{L}).
\end{enumerate}
As a first step towards point (iii) and using some techniques of Mossel, O'Donnel and Oleszkiewicz \cite{MOO},  we establish in Theorem \ref{MOO-thm}
that, if $\mu$ is the law of a sequence of multilinear polynomials with low influences, bounded degree and unit variance, then it necessarily admits a density with respect to Lebesgue measure.

\bigskip

Our second example is concerned with sequences belonging to a fixed order Wiener chaos $\mathcal{H}_k$ (with $k\geq 2$)
and when we have convergence in $L^2(\Omega)$.
More precisely, let $\{F_n\}$ be a sequence of the form $F_n=I_k(f_n)$
(with $I_k$ the $k$th multiple Wiener-It\^o integral)
and assume that it
converges in $L^2(\Omega)$ towards a random variable $F_{\infty}=I_k(f_\infty)$.
Assume moreover that $E[F_\infty^2]>0$.
Then, there exists a constant $c>0$ such that, for any $n\geq 1$,
\begin{equation}\label{ineq2}
d_{TV}(F_n,F_\infty)\leq c\, \|f_n-f_\infty\|^{\frac{1}{2k}}.
\end{equation}
Actually, the inequality (\ref{ineq2}) is not new. It was shown in 1987 by Davydov and Martynova in \cite{DM} (with the better factor $\frac1k$
instead of $\frac1{2k}$). However, it is a pity that \cite{DM} contains only a sketch of the proof of (\ref{ineq2}). Since it is not clear (at least for us!) how to complete the missing details, we believe that
our proof may be of interest as it is fully self-contained.
Moreover, we are hopeful that our approach could be used in the multivariate framework as well, which
would solve an open problem (see indeed \cite{breton} and comments therein). Once again, we postpone this possible extension in a
subsequent paper.

\bigskip

Finally, we develop a third example. It arises when one seeks for a multidimensional counterpart of
Theorem \ref{N-P}, that is, when one wants to prove that one can replace for free the convergence in distribution in
the statement of the Peccati-Tudor theorem (\cite[Theorem 6.2.3]{NouPecBook}) by a convergence in total variation.
We prove, without relying to Stein's method but in the same spirit as in the famous proof of the H\"ormander theorem by Paul Malliavin
\cite{MalliavinHormander}, that if a sequence of vectors of multiple Wiener-It\^o
integrals converges in law to a
Gaussian vector having a  non-degenerate covariance matrix, then it necessarily converges in total variation. This result solves, in the multidimensional framework, a  problem left open after the discovery of Theorem \ref{N-P}.

\bigskip

Our paper contains results closely connected to those of the paper \cite{HuLuNu} by Hu, Lu and Nualart.
The investigations were done independently and at about the same time.
In \cite{HuLuNu}, the authors focus on the convergence of random vectors $\{F_n\}$ which are functionals of Gaussian processes to a normal $\mathcal{N}(0,I_d)$.
More specifically, they work under a negative moment condition (in the spirit of our Theorem \ref{PMa} and whose validity
may be sometimes difficult to check in concrete situations)
which enables them to show that the density of $F_n$ (as well as its first derivatives) converges to the Gaussian density.
Applications to sequences of random variables in the second Wiener chaos is then discussed.
It is worth mentioning that the philosophy of our paper is a bit different. We are indeed interested in exhibiting instances
 for which, {\it without} further assumptions, the convergence in law (to a random variable which is {\it not necessarily} Gaussian)
 turns out to be equivalent to the convergence in total variation\footnote{When we are dealing with sequences of random variables
 that have a law which is absolutely continuous with respect to the Lebesgue measure,
 which is going to be always the case in our paper, it is worthwhile noting that the convergence in total variation is actually equivalent to the $L^1$-convergence of densities.}.

\bigskip

The rest of the paper is organized as follows. In Section \ref{s:preliminaries}, we first recall some useful facts about multiple Wiener-Itô integrals and Malliavin calculus.
We then prove inequality (\ref{ineq}) in Section \ref{s:proofofineq}.
The proof of (\ref{ineq2}) is done in Section \ref{s:davyma}.
Finally, our extension of the Peccati-Tudor Theorem is given in Section \ref{s:pt}.

\section{Preliminaries}\label{s:preliminaries}
This section contains the elements of Gaussian analysis and Malliavin calculus that are used throughout this paper. See the
monographs \cite{NouPecBook, nualartbook} for further details.

\subsection{Isonormal processes and multiple Wiener-It\^o integrals}\label{ss:isonormal}

Let $\EuFrak H$ be a real separable Hilbert space. For any $k\geq 1$, we write $\EuFrak H^{\otimes k}$ and $\EuFrak H^{\odot k}$ to
indicate, respectively, the $k$th tensor power and the $k$th symmetric tensor power of $\EuFrak H$; we also set by convention
$\EuFrak H^{\otimes 0} = \EuFrak H^{\odot 0} =\R$. When $\HH = L^2(A,\mathcal{A}, \mu) =:L^2(\mu)$, where $\mu$ is a $\sigma$-finite
and non-atomic measure on the measurable space $(A,\mathcal{A})$, then $\EuFrak H^{\otimes k} = L^2(A^k,\mathcal{A}^k,\mu^k)=:L^2(\mu^k)$, and $\EuFrak H^{\odot k} = L_s^2(A^k,\mathcal{A}^k,\mu^k) := L_s^2(\mu^k)$, where $L_s^2(\mu^k)$ stands for the subspace of $L^2(\mu^k)$ composed of those functions that are $\mu^k$-almost everywhere symmetric. We denote by $X=\{X(h) : h\in \EuFrak H\}$
an {\it isonormal Gaussian process} over
$\EuFrak H$.
This means that $X$ is a centered Gaussian family, defined on some probability space $(\Omega ,\mathcal{F},P)$, with a covariance structure given by the relation
$E\left[ X(h)X(g)\right] =\langle h,g\rangle _{\EuFrak H}$. We also assume that $\mathcal{F}=\sigma(X)$, that is,
$\mathcal{F}$ is generated by $X$.

For every $k\geq 1$, the symbol $\mathcal{H}_{k}$ stands for the $k$th {\it Wiener chaos} of $X$,
defined as the closed linear subspace of $L^2(\Omega ,\mathcal{F},P) =:L^2(\Omega) $
generated by the family $\{H_{k}(X(h)) : h\in \EuFrak H,\left\|
h\right\| _{\EuFrak H}=1\}$, where $H_{k}$ is the $k$th Hermite polynomial
given by
\begin{equation}\label{hq}
H_k(x) = (-1)^k e^{\frac{x^2}{2}}
 \frac{d^k}{dx^k} \big( e^{-\frac{x^2}{2}} \big).
\end{equation}
We write by convention $\mathcal{H}_{0} = \mathbb{R}$. For
any $k\geq 1$, the mapping $I_{k}(h^{\otimes k})=H_{k}(X(h))$ can be extended to a
linear isometry between the symmetric tensor product $\EuFrak H^{\odot k}$
(equipped with the modified norm $\sqrt{k!}\left\| \cdot \right\| _{\EuFrak H^{\otimes k}}$)
and the $k$th Wiener chaos $\mathcal{H}_{k}$. For $k=0$, we write $I_{0}(c)=c$, $c\in\mathbb{R}$. A crucial fact is that, when $\HH = L^2(\mu)$, for every $f\in\EuFrak H^{\odot k} = L_s^2(\mu^k)$ the random variable $I_k(f)$ coincides with the $k$-fold multiple Wiener-It\^o stochastic integral of $f$ with respect to the centered Gaussian measure (with control $\mu$) canonically generated by $X$ (see \cite[Section 1.1.2]{nualartbook}).

It is
well-known that $L^2(\Omega)$
can be decomposed into the infinite orthogonal sum of the spaces $\mathcal{H}_{k}$. It follows that any square-integrable random variable
$F\in L^2(\Omega)$ admits the following {\it Wiener-It\^{o} chaotic expansion}
\begin{equation}
F=\sum_{k=0}^{\infty }I_{k}(f_{k}),  \label{E}
\end{equation}
where $f_{0}=E[F]$, and the $f_{k}\in \EuFrak H^{\odot k}$, $k\geq 1$, are
uniquely determined by $F$. For every $k\geq 0$, we denote by $J_{k}$ the
orthogonal projection operator on the $k$th Wiener chaos. In particular, if
$F\in L^2(\Omega)$ is as in (\ref{E}), then
$J_{k}F=I_{k}(f_{k})$ for every $k\geq 0$.

Let $\{e_{i},\,i\geq 1\}$ be a complete orthonormal system in $\EuFrak H$.
Given $f\in \EuFrak H^{\odot k}$ and $g\in \EuFrak H^{\odot l}$, for every
$r=0,\ldots ,k\wedge l$, the \textit{contraction} of $f$ and $g$ of order $r$
is the element of $\EuFrak H^{\otimes (k+l-2r)}$ defined by
\begin{equation}
f\otimes _{r}g=\sum_{i_{1},\ldots ,i_{r}=1}^{\infty }\langle
f,e_{i_{1}}\otimes \ldots \otimes e_{i_{r}}\rangle _{\EuFrak H^{\otimes
r}}\otimes \langle g,e_{i_{1}}\otimes \ldots \otimes e_{i_{r}}
\rangle_{\EuFrak H^{\otimes r}}.  \label{v2}
\end{equation}
Notice that the definition of $f\otimes_r g$ does not depend
on the particular choice of $\{e_i,\,i\geq 1\}$, and that
$f\otimes _{r}g$ is not necessarily symmetric; we denote its
symmetrization by $f\widetilde{\otimes }_{r}g\in \EuFrak H^{\odot (k+l-2r)}$.
Moreover, $f\otimes _{0}g=f\otimes g$ equals the tensor product of $f$ and
$g$ while, for $k=l$, $f\otimes_{k}g=\langle f,g\rangle _{\EuFrak H^{\otimes k}}$.
When $\HH = L^2(A,\mathcal{A},\mu)$ and $r=1,...,k\wedge l$, the contraction $f\otimes _{r}g$ is the element of $L^2(\mu^{k+l-2r})$ given by
\begin{eqnarray}\label{e:contraction}
&& f\otimes _{r}g (x_1,...,x_{k+l-2r})\\
&& = \int_{A^r} f(x_1,...,x_{k-r},a_1,...,a_r)g(x_{k-r+1},...,x_{k+l-2r},a_1,...,a_r)d\mu(a_1)...d\mu(a_r). \notag
\end{eqnarray}
It can also be shown that the following {\sl product formula} holds: if $f\in \EuFrak
H^{\odot k}$ and $g\in \EuFrak
H^{\odot l}$, then
\begin{eqnarray}\label{multiplication}
I_k(f) I_l(g) = \sum_{r=0}^{k \wedge l} r! {k \choose r}{ l \choose r} I_{k+l-2r} (f\widetilde{\otimes}_{r}g).
\end{eqnarray}
\smallskip
Finally, we state a very useful property of Wiener chaos (see \cite{E.N} or \cite[Corollary 2.8.14]{NouPecBook}), which is going to be
used several times in the sequel (notably
in the proofs of Lemmas \ref{lemme algèbre} and \ref{paley}).
\begin{theo}[Hypercontractivity]\label{hyper}
Let $F\in \mathcal{H}_k$ with $k\geq 1$. Then, for all $r>1$,
\[
\esp\left[|F|^r\right]^{1/r}\leq
(r-1)^{k/2} E[F^2]^{1/2}.
\]
\end{theo}

\subsection{Malliavin calculus}
We now introduce some basic elements of the Malliavin calculus with respect
to the isonormal Gaussian process $X=\{X(h),\,h\in\HH\}$. Let $\mathcal{S}$
be the set of all
cylindrical random variables of
the form
\begin{equation}
F=g\left( X(\phi _{1}),\ldots ,X(\phi _{n})\right) ,  \label{v3}
\end{equation}
where $n\geq 1$, $g:\mathbb{R}^{n}\rightarrow \mathbb{R}$ is an infinitely
differentiable function such that its partial derivatives have polynomial growth, and $\phi _{i}\in \EuFrak H$,
$i=1,\ldots,n$.
The {\it Malliavin derivative}  of $F$ with respect to $X$ is the element of
$L^2(\Omega ,\EuFrak H)$ defined as
\begin{equation*}
DF\;=\;\sum_{i=1}^{n}\frac{\partial g}{\partial x_{i}}\left( X(\phi
_{1}),\ldots ,X(\phi _{n})\right) \phi _{i}.
\end{equation*}
In particular, $DX(h)=h$ for every $h\in \EuFrak H$. By iteration, one can
define the $m$th derivative $D^{m}F$, which is an element of $L^2(\Omega ,\EuFrak H^{\odot m})$
for every $m\geq 2$.
For $m\geq 1$ and $p\geq 1$, ${\mathbb{D}}^{m,p}$ denotes the closure of
$\mathcal{S}$ with respect to the norm $\Vert \cdot \Vert _{m,p}$, defined by
the relation
\begin{equation*}
\Vert F\Vert _{m,p}^{p}\;=\;\esp\left[ |F|^{p}\right] +\sum_{i=1}^{m}\esp\left[
\Vert D^{i}F\Vert _{\EuFrak H^{\otimes i}}^{p}\right].
\end{equation*}
We often use the notation $\mathbb{D}^{\infty} := \bigcap_{m\geq 1}
\bigcap_{p\geq 1}\mathbb{D}^{m,p}$.

\begin{rem}\label{r:density}{\rm Any random variable $Y$ that is a finite linear combination of multiple Wiener-It\^o integrals is an element of $\mathbb{D}^\infty$.
Moreover, if $Y\neq 0$, then the law of $Y$ admits a density with respect to the Lebesgue measure -- see \cite{Shigekawa}
or \cite[Theorem 2.10.1]{NouPecBook}.
}
\end{rem}

The Malliavin derivative $D$ obeys the following \textsl{chain rule}. If
$\varphi :\mathbb{R}^{n}\rightarrow \mathbb{R}$ is continuously
differentiable with bounded partial derivatives and if $F=(F_{1},\ldots
,F_{n})$ is a vector of elements of ${\mathbb{D}}^{1,2}$, then $\varphi
(F)\in {\mathbb{D}}^{1,2}$ and
\begin{equation}\label{e:chainrule}
D\,\varphi (F)=\sum_{i=1}^{n}\frac{\partial \varphi }{\partial x_{i}}(F)DF_{i}.
\end{equation}

\begin{rem}\label{r:polg}{\rm
By approximation, it is easily checked that equation (\ref{e:chainrule}) continues to hold in the following two cases:
(i) $F_i \in \mathbb{D}^\infty$ and $\varphi$ has continuous partial derivatives with at most polynomial growth, and
(ii) $F_i\in \mathbb{D}^{1,2}$ has an absolutely continuous distribution and $\varphi$ is Lipschitz continuous.
}
\end{rem}

Note also that a random variable $F$ in $L^2(\Omega)$ is in ${\mathbb{D}}^{1,2}$ if and only if
$\sum_{k=1}^{\infty }k\|J_kF\|^2_{L^2(\Omega)}<\infty$
and, in this case, $\esp\left[ \Vert DF\Vert _{\EuFrak H}^{2}\right]
=\sum_{k=1}^{\infty }k\|J_kF\|^2_{L^2(\Omega)}$. If $\EuFrak H=
L^{2}(A,\mathcal{A},\mu )$ (with $\mu $ non-atomic), then the
derivative of a random variable $F$ in $L^2(\Omega)$ can be identified with
the element of $L^2(A \times \Omega )$ given by
\begin{equation}
D_{x}F=\sum_{k=1}^{\infty }kI_{k-1}\left( f_{k}(\cdot ,x)\right) ,\quad x \in A.  \label{dtf}
\end{equation}

We denote by $\delta $ the adjoint of the operator $D$, also called the
\textit{divergence operator}. A random element $u\in L^2(\Omega ,\EuFrak H)$
belongs to the domain of $\delta $, noted $\mathrm{Dom}\,\delta $, if and
only if it verifies
$|\esp\langle DF,u\rangle _{\EuFrak H}|\leq c_{u}\,\Vert F\Vert _{L^2(\Omega)}$
for any $F\in \mathbb{D}^{1,2}$, where $c_{u}$ is a constant depending only
on $u$. If $u\in \mathrm{Dom}\,\delta $, then the random variable $\delta (u)$
is defined by the duality relationship (customarily called \textit{integration by parts
formula})
\begin{equation}
\esp[F\delta (u)]=\esp[\langle DF,u\rangle _{\EuFrak H}],  \label{ipp0}
\end{equation}
which holds for every $F\in {\mathbb{D}}^{1,2}$.
More generally, if $F\in \mathbb{D}^{1,2}$ and $u\in {\rm Dom}\,\delta$ are such that the three expectations
$\esp\big[F^2\|u\|^2_\HH]$, $\esp\big[F^2\delta(u)^2\big]$ and
$\esp\big[\langle DF,u\rangle_\HH^2\big]$ are finite, then
 $Fu\in{\rm Dom}\,\delta$ and
\begin{equation}\label{etoile}
\delta(Fu)=F\delta(u)- \langle DF,u\rangle_\HH.
\end{equation}

The operator $L$, defined as
$L=\sum_{k=0}^{\infty }-kJ_{k}$,
is the {\it infinitesimal generator of the Ornstein-Uhlenbeck
semigroup}.
The domain of $L$ is
\begin{equation*}
\mathrm{Dom}L=\{F\in L^2(\Omega ):\sum_{k=1}^{\infty }k^{2}\left\|
J_{k}F\right\| _{L^2(\Omega )}^{2}<\infty \}=\mathbb{D}^{2,2}\text{.}
\end{equation*}
There is an important relation between the operators $D$, $\delta $ and $L$.
A random variable $F$ belongs to
$\mathbb{D}^{2,2}$ if and only if $F\in \mathrm{Dom}\left( \delta D\right) $
(i.e. $F\in {\mathbb{D}}^{1,2}$ and $DF\in \mathrm{Dom}\delta $) and, in
this case,
\begin{equation}
\delta DF=-LF.  \label{k1}
\end{equation}
In particular, if $F\in\mathbb{D}^{2,2}$ and $H,G\in\mathbb{D}^{1,2}$ are such that $HG\in\mathbb{D}^{1,2}$, then
\begin{eqnarray}\notag
-\esp[HG\,LF]&=&\esp[HG\,\delta DF]=\esp[\langle D(HG),DF\rangle_\HH]=\esp[H\langle DG,DF\rangle_\HH]+\esp[G\langle DH,DF\rangle_\HH].\\
\label{voiture}
\end{eqnarray}

\subsection{A useful result}

In this section, we state and prove the following lemma, which will be used several times
in the sequel.
\begin{lem}\label{paley}
Fix $p\geq 2$, and let $\{F_n\}$ be a sequence of non-zero random variables belonging to the finite sum of chaoses $\bigoplus_{k=0}^p\mathcal{H}_k$.
Assume that $F_n$ converges in distribution as $n\to\infty$.
Then $\sup_{n\geq 1}\esp[|F_n|^r]<\infty$ for all $r\geq 1$.
\end{lem}
{\it Proof}. Let $Z$ be a positive random variable such that $\esp[Z]=1$.
Consider the decomposition $Z=Z{\bf 1}_{\{Z\geq 1/2\}}+Z{\bf 1}_{\{Z<1/2\}}$ and
take the expectation. One deduces, using Cauchy-Schwarz, that
\[
1\leq \sqrt{\esp[Z^2]}\sqrt{P(Z\geq 1/2)}+\frac12,
\]
that is,
\begin{equation}\label{papalele}
\esp[Z^2]\,P(Z\geq 1/2)\geq \frac14.
\end{equation}
On the other hand, Theorem \ref{hyper} implies the existence of $c_p>0$ (a constant depending only on $p$)
such that
$\esp\big[F_n^4\big]\leq c_p\,\esp\big[F_n^2\big]^2$ for all $n\geq 1$.
Combining this latter fact with (\ref{papalele})
yields, with $Z=F_n^2/\esp[F_n^2]$,
\begin{equation}\label{Paley}
P\big( F_n^2 \geq \frac12 \esp [F_n^2] \big) \geq \frac1{4c_p}.
\end{equation}
The sequence $\{F_n\}_{n\geq 1}$ converging in distribution, it is tight and one can choose $M>0$
 large enough so that $P(F_n^2 > M) < \frac{1}{4c_p}$
for all $n\geq 1$.
By applying \eqref{Paley}, one obtains that
\begin{align*}
P\big(F_n^2 \geq M\big) < \frac{1}{4c_p} \leq P\big( F_n^2 \geq \frac12\, \esp [F_n^2] \big),
\end{align*}
from which one deduces immediately that $\sup_{n\geq 1} \esp [F_n^2] \leq 2M<\infty$.
The desired conclusion follows from Theorem \ref{hyper}.
\qed

\subsection{Carbery-Wright inequality}
The proof of (\ref{ineq2}) shall rely on the following nice inequality due to Carbery and Wright \cite{CW}.
We state it in the case of standard Gaussian random variables only. But its statement is actually more general, as it works under
a log-concave density assumption.

\begin{thm}[Carbery-Wright]\label{cw-thm}
There exists an absolute constant $c>0$ such that, for all polynomial $Q:\R^n\to \R$ of degree at most $d$,
all independent random variables $X_1,\ldots,X_n\sim \mathcal{N}(0,1)$ and all $\alpha>0$,
\begin{equation}\label{cw-ineq}
E[Q(X_1,\ldots,X_n)^2]^{\frac1{2d}}\,\,
P(|Q(X_1,\ldots,X_n)|\leq \alpha)
\leq c\,d\,\alpha^{\frac1d}.
\end{equation}
\end{thm}
{\it Proof}.
See \cite[Theorem 8]{CW}.
\qed

\bigskip

The power of $\alpha$ in the inequality (\ref{cw-ineq}) is sharp. To see it, it suffices to consider the case where $n=1$ and $Q(x)=x^d$; we then
have
\[
P(|X_1|^d\leq \alpha)=P(|X_1|\leq \alpha^{1/d})\sim_{\alpha\to 0^+}\sqrt{\frac2\pi}\,\alpha^{\frac1d}.
\]

\section{An asymptotic version of a theorem by Shigekawa}\label{s:proofofineq}

Our first result, which may be seen as an asymptotic version of Shigekawa \cite{Shigekawa}, reads
as follows.

\begin{theo}\label{Loi donne VT}
Fix $p\geq 2$, and let $\{F_n\}$ be a sequence of random variables belonging to the finite sum of chaoses $\bigoplus_{k=0}^p\mathcal{H}_k$.
 Assume that $F_n$ converges in distribution towards $F_\infty$ satisfying ${\rm Var}(F_\infty)>0$. Then,
 the following three assertions hold true:
 \begin{enumerate}
 \item
 the sequence $\{F_n\}$ is uniformly bounded in all the $L^r(\Omega)$: that is,  $\sup_{n\geq 1}\esp[|F_n|^r]<\infty$ for all $r\geq 1$;
 \item
 there exists
 $c>0$ such that, for all $n\geq 1$,
\begin{equation}\label{Final}
d_{TV}(F_n,F_\infty)\leq c\, d_{FM}(F_n,F_\infty)^{\frac{1}{2p+1}}.
\end{equation}
In particular, $F_n$ converges in total variation towards $F_\infty$;
 \item
 the law of $F_\infty$ is absolutely continuous with respect to the Lebesgue measure.
 \end{enumerate}
\end{theo}
{\it Proof}.
The first point comes directly from Lemma \ref{paley}.
The rest of the proof is divided into four steps.
Throughout the proof, the letter $c$ stands for a non-negative constant {\it independent} of $n$
(but which may depend on $p$, $\{F_n\}$ or $F_\infty$) and whose value may change
from line to line.
\bigskip

{\it First step}.
We claim that there exists $c>0$ such that, for all $n\geq 1$:
\begin{equation}\label{claim1bis}
P(\|D F_n\|_{\HH}\leq \lambda)\leq c\, \frac{\lambda^{\frac1{p-1}}}{\var(F_n)^{\frac1{2p-2}}}.
\end{equation}
Indeed, let $f_{k,n}$ be the elements of $\HH^{\odot k}$ such that $F_n=E[F_n]+\sum_{k=1}^p I_k(f_{k,n})$. Using the product formula (\ref{multiplication}), we can write:
\begin{eqnarray*}
\|D F_n\|^2_{\HH}
&=&
\sum_{k,l=1}^p kl\langle I_{k-1}(f_{k,n}),I_{l-1}(f_{l,n})\rangle_\HH\\
&=&
\sum_{k,l=1}^p kl\sum_{r=1}^{k\wedge l}(r-1)!\binom{k-1}{r-1}\binom{l-1}{r-1}I_{k+l-2r}(f_{k,n}\tilde{\otimes}_r f_{l,n}).
\end{eqnarray*}
Now, let $\{e_i\}_{i\geq 1}$ be an orthonormal family of $\HH$ and decompose
\[
f_{k,n}\tilde{\otimes}_r f_{l,n}=\displaystyle{\sum_{m_1,m_2,\cdots,m_{k+l-2r}=1}^{\infty}\alpha_{m_1,\cdots,m_{k+l-2r},n}\,e_{m_1}\otimes\cdots\otimes e_{m_{k+l-2r}}}.
\]
Also, set
\[
\begin{array}{ccc}
g_{k,l,r,n,s}&=&\displaystyle{\sum_{m_1,m_2,\cdots,m_{k+l-2r}=1}^{s}\alpha_{m_1,\cdots,m_{k+l-2r},n}\,e_{m_1}\otimes\cdots\otimes e_{m_{k+l-2r}}},\\\\
Y_{s,n}&=&\displaystyle{\sum_{k,l=1}^p kl\sum_{r=1}^{k\wedge l}(r-1)!\binom{k-1}{r-1}\binom{l-1}{r-1}I_{k+l-2r}(g_{k,l,r,n,s})}.
\end{array}
\]
Firstly, it is clear that $g_{k,l,r,n,s}\to f_{k,n}\tilde{\otimes}_r f_{l,n}$ as $s$ tends to infinity in
$\HH^{\otimes (k+l-2r)}$. Hence,
using the isometry property of Wiener-Itô integrals we conclude that
\begin{equation}\label{Y_s tend vers zero}
Y_{s,n}\overset{L^2}{\to}\|DF_n\|^2_\HH\,\mbox{ as $s\to\infty$}.
\end{equation}
We deduce that there exists a strictly increasing sequence $\{s_l\}$ such that $Y_{s_l,n}\to \|DF_n\|^2_\HH$ as $l\to\infty$ almost surely.
Secondly, we deduce from a well-known result by It\^o that, with $k=k_1+\ldots+k_m$,
\[
I_k(e_1^{\otimes k_1}\otimes \ldots\otimes e_m^{\otimes k_m})=\prod_{i=1}^m H_{k_i}\left(X(e_i)\right).
\]
Here, $H_k$ stands for the $k$th Hermite polynomial and has degree $k$, see (\ref{hq}).
Also, one should note that the value of $I_k(e_1^{\otimes k_1}\otimes \ldots\otimes e_m^{\otimes k_m})$
is not modified when one permutes the order of the elements in the tensor product.
Putting these two facts together, we can write
\[
Y_{s,n}=Q_{s,n}\left(X(e_1),\ldots,X(e_s)\right),
\]
for some polynomial $Q_{s,n}$ of degree at most $2p-2$.
Consequently, we deduce from Theorem \ref{cw-thm}
that there exists a constant $c>0$ such that, for any $n\geq 0$ and any $\lambda>0$,
\[
P(\left|Y_{s,n}\right|\leq \lambda^2)
\leq c\, E[Y_{s,n}^2]^{-\frac{1}{4p-4}}\, \lambda^{1/(p-1)}.
\]
Next, we can use Fatou's lemma to deduce that, for any $n\geq 0$ and any $\lambda>0$,
\begin{eqnarray}\notag
P(\|DF_n\|_{\HH}\leq \lambda)&\leq&
P\left(\liminf_{l\to\infty} \{|Y_{s_l,n}|\leq 2\lambda^2 \}\right)\\
&\leq&\liminf_{l\to\infty} P(|Y_{s_l,n}|\leq 2\lambda^2)\leq c\, \esp\left[\|DF_n\|^4_{\HH}\right]^{-1/(4p-4)} \lambda^{1/(p-1)}.\label{intermédiaire}
\end{eqnarray}
Finally, by applying
the Poincar\'e inequality
(that is, $\var(F_n)\leq E[\|DF_n\|^2_\HH]$), we get
\[
\esp\left[\|DF_n\|^4_{\HH}\right]\geq\esp\left[\|DF_n\|^2_{\HH}\right]^2\geq {\rm Var}(F_n)^2,
\]
which, together with (\ref{intermédiaire}), implies the desired conclusion (\ref{claim1bis}).

\bigskip

{\it Second step}.
We claim that: $(i)$ the law of $F_n$ is absolutely continuous with respect to the Lebesgue measure when $n$ is large enough, and that $(ii)$ there exists $c>0$ and $n_0\in\N$ such that, for all $\varepsilon>0$, the following inequality holds:
\begin{equation}\label{claim1}
\sup_{n\geq n_0} \esp\left[\frac{\varepsilon}{\|D F_n\|_{\HH}^2+\varepsilon}\right]\leq c\, \varepsilon^{\frac{1}{2p-1}}.
\end{equation}
The first point is a direct consequence of  Shigekawa \cite{Shigekawa}, but one can also give a direct
proof by using (\ref{claim1bis}). Indeed Point 1 together with the assumption that $F_n$ converges in distribution to $F_\infty$ implies
that $\var(F_n)\to \var(F_\infty)>0$. By letting $\lambda\to 0$ in (\ref{claim1bis}), we deduce that, for $n$ large enough
(so that ${\rm Var}(F_n)>0$), we have $P(\|DF_n\|_\HH=0)=0$. Then, the Bouleau-Hirsch criterion (see, e.g., \cite[Theorem 2.1.3]{nualartbook}) ensures that  the law of $F_n$ is absolutely continuous with respect to the Lebesgue measure.

Now, let us prove (\ref{claim1}).
We deduce from (\ref{claim1bis}) that, for any $\lambda,\e>0$,
\begin{eqnarray*}
\esp\left[\frac{\varepsilon}{\|D F_n\|_{\HH}^2+\varepsilon}\right]&\leq&\esp\left[\frac{\varepsilon}{\|D F_n\|_{\HH}^2+\varepsilon}\,\mathbf{1}_{\{\|D F_n\|_{\HH}>\lambda\}}\right]+P(\|D F_n\|_{\HH}\leq\lambda)\\
&\leq&\frac{\varepsilon}{\lambda^2}+c\, \var(F_n)^{-1/(2p-2)} \lambda^{1/(p-1)}.
\end{eqnarray*}
As we said, we have that $\var(F_n)\to\var(F_\infty)>0$ as
$n\to\infty$.
Therefore, there exists $a>0$ such that $\var(F_n)>a$ for $n$ large enough (say $n\geq n_0$).
We deduce
that there exists $c>0$  such that, for any $\lambda,\e>0$,
\begin{equation}\label{final}
\sup_{n\geq n_0}\esp\left[\frac{\varepsilon}{\|D F_n\|_{\HH}^2+\varepsilon}\right]\leq c \left(\frac{\varepsilon}{\lambda^2}+\lambda^{1/(p-1)}\right).
\end{equation}
Choosing $\lambda=\e^{\frac{p-1}{2p-1}}$ concludes the proof of (\ref{claim1}).

\bigskip

{\it Third step}.
We claim that there exists $c>0$ such that, for all $n,m$ large enough,
 \begin{equation}\label{FactC}
d_{TV}(F_n,F_m)\leq c\, d_{FM}(F_n,F_m)^{\frac{1}{2p+1}}.
\end{equation}
Set $p_{\alpha}(x)=\frac{1}{\alpha\sqrt{2\pi}}e^{-\frac{x^2}{2\alpha^2}}$, $x\in\R$, $0<\alpha\leq 1$.
Let $A$ be a bounded Borel set.
It is easily checked that
\begin{equation}\label{fc1}
\|{\bf 1}_{A}*p_\alpha\|_\infty\leq \|{\bf 1}_{A}\|_\infty \|p_\alpha\|_1=1\leq \frac1\alpha
\end{equation}
and, since $p_\alpha'(x)=-\frac{x}{\alpha^2}p_\alpha(x)$, that
\begin{eqnarray}
\|({\bf 1}_{A}*p_\alpha)'\|_\infty&=&
\|{\bf 1}_{A}*p'_\alpha\|_\infty \notag
=\frac{1}{\alpha^2} \sup_{x\in\R}\left| \int_\R {\bf 1}_{A}(x-y)\,y\,p_\alpha(y)dy\right|\\
&\leq&
\frac{1}{\alpha^2} \int_\R |y|p_\alpha(y)dy=\frac{1}{\alpha}\sqrt{\frac2\pi}\leq\frac1\alpha.\label{Fact-convo}
\end{eqnarray}
Let $n,m$ be large integers. Using (\ref{claim1}), (\ref{fc1}), (\ref{Fact-convo})
and that $F_n$ has a density when $n$ is large enough (Step 2 $(i)$), we can write
\begin{eqnarray*}
&&\big|P(F_n\in A)-P(F_m\in A)\big|\\
&\leq&\big|\esp\left[{\bf 1}_A*p_\alpha(F_n)-{\bf 1}_A*p_\alpha(F_m)\right]\big|\\
&&
+\left|\esp\left[\left({\bf 1}_A(F_n)-{\bf 1}_A*p_\alpha(F_n)\right)\left(\frac{\|D F_n\|_{\HH}^2}{\|D F_n\|_{\HH}^2+\varepsilon}+\frac{\e}{\|D F_n\|_{\HH}^2+\varepsilon}\right)\right]\right|\\
&&
+\left|\esp\left[\left({\bf 1}_A(F_m)-{\bf 1}_A*p_\alpha(F_m)\right)\left(\frac{\|D F_m\|_{\HH}^2}{\|D F_m\|_{\HH}^2+\varepsilon}+\frac{\e}{\|D F_m\|_{\HH}^2+\varepsilon}\right)\right]\right|\\
&\leq & \frac{1}{\alpha}d_{FM}(F_n,F_m)+2\esp\left[\frac{\varepsilon}{\|D F_n\|_{\HH}^2+\varepsilon}\right]
+2\esp\left[\frac{\varepsilon}{\|D F_m\|_{\HH}^2+\varepsilon}\right]\\
&&
+\left|\esp\left[\left({\bf 1}_A(F_n)-{\bf 1}_A*p_\alpha(F_n)\right)\frac{\|D F_n\|_{\HH}^2}{\|D F_n\|_{\HH}^2+\varepsilon}\right]\right|
+\left|\esp\left[\left({\bf 1}_A(F_m)-{\bf 1}_A*p_\alpha(F_m)\right)\frac{\|D F_m\|_{\HH}^2}{\|D F_m\|_{\HH}^2+\varepsilon}\right]\right|
\\
&\leq & \frac{1}{\alpha}d_{FM}(F_n,F_m)+c\,\e^{1/(2p-1)}
+2\sup_{n\geq n_0}\left|\esp\left[\left({\bf 1}_A(F_n)-{\bf 1}_A*p_\alpha(F_n)\right)\frac{\|D F_n\|_{\HH}^2}{\|D F_n\|_{\HH}^2+\varepsilon}\right]\right|.
\end{eqnarray*}
Now, set $\Psi(x)=\int_{-\infty}^{x} {\bf 1}_A(s) ds$ and let us integrate by parts through (\ref{voiture}). We get
\begin{eqnarray}
\notag
&&\left|\esp\left[\left({\bf 1}_A(F_n)-{\bf 1}_A*p_\alpha(F_n)\right)\frac{\|D F_n\|_{\HH}^2}{\|D F_n\|_{\HH}^2+\varepsilon}\right]\right|\\
&=&\left|\esp\left[\frac{1}{\|D F_n\|_{\HH}^2+\varepsilon}\,\langle D(\Psi(F_n)-\Psi*p_{\alpha}(F_n)),DF_n\rangle_\HH\right]\right|\notag\\
\notag
&=&\left|\esp\left[\left(\Psi(F_n)-\Psi*p_{\alpha}(F_n)\right)\left(\left\langle DF_n,D\left(\frac{1}{\|D F_n\|_{\HH}^2+\varepsilon}\right)\right\rangle_\HH+\frac{LF_n}{\|D F_n\|_{\HH}^2+\varepsilon}\right)\right]\right|\\\notag
&=&\left|\esp\left[\left(\Psi(F_n)-\Psi*p_{\alpha}(F_n)\right)\left(-\frac{2\langle D^2F_n,DF_n\otimes DF_n\rangle_{\HH^{\otimes 2}}}{(\|D F_n\|_{\HH}^2+\varepsilon)^2}+\frac{LF_n}{\|D F_n\|_{\HH}^2+\varepsilon}\right)\right]\right|\\
&\leq&\frac1\e\,\esp\left[\left|\Psi(F_n)-\Psi*p_{\alpha}(F_n)\right|\big(2\|D^2F_n\|_{\HH^{\otimes 2}}+\big|LF_n\big|\big)\right].\label{FactA}
\end{eqnarray}
On the other hand, we have
\begin{eqnarray}
\notag
\left|\Psi(x)-\Psi*p_{\alpha}(x)\right|&=&\left|\int_\R p_\alpha(y)\left(\int_{-\infty}^x \left({\bf 1}_A(z)-{\bf 1}_A(z-y)\right)dz\right)dy\right|\\\notag
&\leq&\int_\R p_\alpha(y)\left|\int_{-\infty}^x {\bf 1}_A(z)dz-\int_{-\infty}^x{\bf 1}_A(z-y)dz\right|dy\\
&\leq&\int_\R p_\alpha(y)\left|\int_{x-y}^x {\bf 1}_A(z)dz\right| dy
\leq \int_\R p_{\alpha}(y)\left|y\right|dy \leq \sqrt{\frac{2}{\pi}}\alpha.\label{FactB}
\end{eqnarray}
Moreover, $F_n$ is bounded in $L^2(\Omega)$ (see indeed Point 1) and $D^2 F_n=\sum_{k=2}^p k(k-1)I_{k-2}(f_{k,n})$. We deduce that $\sup_{n\geq 1}E[|LF_n|]<\infty$ (since $F_n\in\bigoplus_{k=0}^p \mathcal{H}_k$) and
\[\sup_n\esp\left[\|D^2F_n\|_{\HH^{\otimes 2}}\right]<\infty,\]
implying in turn, thanks to (\ref{FactB}), that
\begin{eqnarray*}
\sup_{n\geq n_0}\left|\esp\left[\left({\bf 1}_A(F_n)-{\bf 1}_A*p_\alpha(F_n)\right)\frac{\|D F_n\|_{\HH}^2}{\|D F_n\|_{\HH}^2+\varepsilon}\right]\right|\leq c\,\frac{\alpha}{\varepsilon}.
\end{eqnarray*}
Thus, there exists $c>0$ and $n_0\in\N$ such that, for any $n,m\geq n_0$, any $0<\alpha\leq 1$ and any $\e>0$,
\begin{eqnarray*}
d_{TV}(F_n,F_m)
\leq c\left(\frac{1}{\alpha}d_{FM}(F_n,F_m)+\frac{\alpha}{\varepsilon}+\varepsilon^{\frac{1}{2p-1}}\right).
\end{eqnarray*}
Choosing $\alpha=\big(\frac12 d_{FM}(F_n,F_m)\big)^{\frac{2p}{2p+1}}$ (observe that $\alpha\leq 1$) and $\e=d_{FM}(F_n,F_m)^{\frac{2p-1}{2p+1}}$ leads to
our claim
(\ref{FactC}).

\bigskip

{\it Fourth and final step}.
Since the Fortet-Mourier distance $d_{FM}$ metrizes the convergence in distribution (see, e.g., \cite[Theorem 11.3.3]{Dudley book}), our assumption ensures that $d_{FM}(F_n,F_m)\to 0$ as $n,m\to\infty$.
Thanks to (\ref{FactC}), we conclude that $\bP_{F_n}$ is a Cauchy sequence for the total variation distance. But the space of bounded measures is complete for the total variation distance, so $\bP_{F_n}$ must converge towards $\bP_{F_\infty}$ in the total variation distance. Letting $m\to\infty$ in (\ref{FactC}) yields the desired inequality (\ref{Final}). The proof of point 2 is done.

Let $A$ be a Borel set of Lebesgue measure zero.
By Step 2 $(i)$, we have $P(F_n\in A)=0$ when $n$ is large enough.
Since $d_{TV}(F_n,F_\infty)\to 0$ as $n\to\infty$, we deduce that $P(F_\infty\in A)=0$,
proving that the law of $F_\infty$ is absolutely continuous with respect to the Lebesgue measure by
the Radon-Nikodym theorem. The proof of point 3 is done.
\qed

\bigskip

Let us give an application  of Theorem \ref{Loi donne VT} to the study of the absolute continuity of laws which are limits of multilinear polynomials with low influences and bounded degree. We use techniques from  Mossel, O'Donnel and Oleszkiewicz \cite{MOO}.

\begin{thm}\label{MOO-thm}
Let $p\geq 1$ be an integer and let $X_1$, $X_2$, $\ldots$ be independent random variables. Assume further that $E[X_k]=0$ and $E[X_k^2]=1$ for all $k$ and that there exists $\e>0$ such that $\sup_k E|X_k|^{2+\e}<\infty$.
For any $m\geq 1$, let $n=n(m)$ et let $Q_m\in\R[x_1,\ldots,x_{n(m)}$ be
a real polynomial of the form
\[
Q_m(x_1,\ldots,x_{n(m)}) = \sum_{S\subset\{1,\ldots,n(m)\}}c_{S,m}\,\prod_{i\in S}x_i,
\]
with $\sum_{|S|>0}c_{S,m}^2=1$. Suppose moreover that 
the contribution of each $x_i$ to $Q_m(x_1,\ldots,x_{n(m)})$ is uniformly negligible, that is,
\begin{equation}\label{lowinfluence}
\lim_{m\to\infty}\sup_{1\leq k\leq n(m)} \sum_{S:k\in S} c_{S,m}^2 =0,
\end{equation}
and that the degree of $Q_m$ is at most $p$, that is,
\begin{equation}\label{boundeddegree}
\max_{S:\,c_{s,m}\neq 0}|S|\leq p.
\end{equation}
Finally, let $F$ be a limit in law of $Q_m(X_1,\ldots,X_{n(m)})$ (possibly
through a subsequence only) as $m\to\infty$. Then the law of $F$ has a density
with respect to the Lebesgue measure.
\end{thm}
{\it Proof}. Using \cite[Theorem 2.2]{MOO} and because of (\ref{lowinfluence}) and (\ref{boundeddegree}), we deduce that
$F$ is also a limit in law of $Q_m(G_1,\ldots,G_{n(m)})$, where
the $G_i$'s are independent $N(0,1)$ random variables.
Moreover, because of (\ref{boundeddegree}), it is straightforward that $Q_m(G_1,\ldots,G_{n(m)})$ may
be realized as an element belonging to the sum of the $p$ first Wiener chaoses. Also, due to $\sum_{|S|>0}c_{S,m}^2=1$, we have that
the variance of $Q_m(G_1,\ldots,G_{n(m)})$ is 1. Therefore, the desired conclusion is now a direct consequence of Theorem \ref{Loi donne VT} .\qed

\bigskip

\begin{rem}
{\rm
One cannot remove the assumption (\ref{lowinfluence}) in the previous theorem. Indeed, without this assumption, it is straightforward to construct easy counterexamples to the conclusion of Theorem \ref{Loi donne VT}. For instance, it is clear that the conclusion is not reached if one considers $Q_m(x_1,\ldots,x_{n(m)})=x_1$ together with a {\it discrete} random
variable $X_1$.
}
\end{rem}

\section{Continuity of the law of $I_k(f)$ with respect to $f$}\label{s:davyma}

In this section, we are mainly interested in the continuity of the law of $I_k(f)$ with respect to its kernel $f$.
Our first theorem is a result going in the same direction. It exhibits a sufficient condition that allows one to pass from a convergence in law to a convergence in total variation.

\begin{thm}
Let $\{F_n\}_{n\geq 1}$ be a sequence of $\mathbb{D}^{1,2}$ satisfying
\begin{enumerate}
\item[(i)] $\frac{DF_n}{\|DF_n\|_\HH^2}\in {\rm dom}\delta$ for any
$n\geq 1$;
\item[(ii)] $C:=\sup_{n\geq 1}E\left|\delta\left( \frac{DF_n}{\|DF_n\|_\HH^2}\right)\right|<\infty$;
\item[(iii)] $F_n\overset{\rm law}{\to}F_\infty$ as $n\to\infty$.
\end{enumerate}
Then
\[
d_{TV}(F_n,F_\infty)\leq \big(\sqrt{2}+\frac{4C}{\sqrt{\pi}}\big)\sqrt{d_{FM}(F_n,F_\infty)}.
\]
In particular, $F_n$ tends to $F_\infty$ in total variation.
\end{thm}
{\it Proof}.
Let $A$ be a bounded Borel set and
set $p_\alpha(x)=\frac{1}{\sqrt{2\pi}\alpha}e^{-\frac{x^2}{2\alpha^2}}$, $x\in\R$, $0<\alpha\leq 1$.
Since $\int_0^\cdot \big({\bf 1}_{A}(x)-{\bf 1}_{A}*p_\alpha\big)dx$ is Lipschitz
and $F_m$, $F_n$ admit a density,
we have using (\ref{ipp0}) that, for any $n,m\geq 1$,
\begin{eqnarray*}
&&P(F_n\in A) - P(F_m\in A)\\
&=&E[{\bf 1}_{A}*p_\alpha(F_n)] - E[{\bf 1}_{A}*p_\alpha(F_m)]
+ E\left[
\delta
\left(
\frac{DF_n}{\|DF_n\|_\HH^2}
\right)
\int_0^{F_n} \big({\bf 1}_{A}-{\bf 1}_{A}*p_\alpha\big)(x)dx
\right]\\
&&-E\left[
\delta
\left(
\frac{DF_m}{\|DF_m\|_\HH^2}
\right)
\int_0^{F_m} \big({\bf 1}_{A}-{\bf 1}_{A}*p_\alpha\big)(x)dx
\right].
\end{eqnarray*}
Using (\ref{fc1}) and (\ref{Fact-convo}), we can write
\[
\big|E[{\bf 1}_{A}*p_\alpha(F_n)] - E[{\bf 1}_{A}*p_\alpha(F_m)]\big|\leq \frac1\alpha d_{FM}(F_n,F_m).
\]
On the other hand, we have, for any $x\in\R$,
\begin{eqnarray*}
&&\left|\int_0^x \big({\bf 1}_{A}-{\bf 1}_{A}*p_\alpha\big)(v)dv\right|
=\left|\int_0^x du\int_\R dv p_\alpha(v)\big({\bf 1}_{A}(u)-{\bf 1}_{A}(u-v)\big)\right|\\
&=&\left|\int_\R dv\,p_\alpha(v)\int_0^x \big({\bf 1}_{A}(u)-{\bf 1}_{A}(u-v)\big)du\right|\\
&=&\left|\int_\R dv\,p_\alpha(v)\left(\int_{x-v}^x {\bf 1}_{A}(u)du-\int_0^{-v}{\bf 1}_{A}(u)du\right)\right|
\leq 2\int_\R |v|p_\alpha(v)dv =2\sqrt{\frac2\pi}\,\alpha.
\end{eqnarray*}
By putting all these facts together, we get that
\begin{eqnarray*}
\big|P(F_n\in A) - P(F_m\in A)\big|
\leq
\frac1\alpha d_{FM}(F_n,F_m) + 4C\sqrt{\frac2\pi}\alpha.
\end{eqnarray*}
To conclude, it remains to choose $\alpha=\sqrt{\frac12d_{FM}(F_n,F_m)}$ and then to let $m\to\infty$
as in the fourth step of the proof of Theorem \ref{Loi donne VT}.
\qed

\bigskip
In \cite{PM}, Poly and Malicet prove that, if $F_n \overset{\D^{1,2}}{\to} F_\infty$ and $P(\|DF_\infty\|_\HH>0)=1$, then $d_{TV}(F_n,F_\infty)\to 0$. Nevertheless, their proof does not give any idea on the rate of convergence. The following result is a kind of quantitative version of the aforementioned result in \cite{PM}.
\begin{thm}\label{PMa}
Let $\{F_n\}_{n\geq 1}$ be a sequence in $\mathbb{D}^{1,2}$ such that each $F_n$ admits a density.
Let $F_\infty\in\mathbb{D}^{2,4}$ and let $0<\alpha\leq 2$ be such that $\esp\left[\frac{1}{\|DF_\infty\|_\HH^\alpha}\right]<\infty$.
If $F_n\overset{\mathbb{D}^{1,2}}{\to}F_\infty$ then there exists a constant $c>0$ depending only of $F_\infty$
such that, for any $n\geq 1$,
\[d_{TV}(F_n,F_\infty)\leq c \|F_n-F_\infty\|_{\D^{1,2}}^{\frac{\alpha}{\alpha+2}}.\]
\end{thm}
{\it Proof}.
Throughout the proof, the letter $c$ stands for a non-negative constant {\it independent} of $n$
and whose value may change from line to line.
Let $A$ be a bounded Borel set of $\R$. For all $0<\e\leq 1$, one has (using that $F_n$ has a density to perform
the integration by parts, see Remark \ref{r:polg})
\begin{eqnarray}
&&P(F_n\in A)-P(F_\infty\in A)\notag\\
&=&E\left[\frac{\langle D\int_{F_\infty}^{F_n}{\bf 1}_A(x)dx,DF_\infty\rangle_\HH}{\|DF_\infty\|_\HH^2+\e}\right]
+E\left[\big({\bf 1}_A(F_n)-{\bf 1}_A(F_\infty)\big)\frac{\e}{\|DF_\infty\|_\HH^2+\e}\right]\notag\\
&&-E\left[{\bf 1}_A(F_n)\frac{\langle D(F_n-F_\infty),DF_\infty\rangle_\HH}{\|DF_\infty\|_\HH^2+\e}\right]\label{lili1}
\end{eqnarray}
But, see (\ref{etoile}),
\[
\left\langle D\int_{F_\infty}^{F_n}{\bf 1}_A(x)dx,DF_\infty\right\rangle_\HH
=-\delta\left(
DF_\infty \int_{F_\infty}^{F_n}{\bf 1}_A(x)dx
 \right)
 +LF_\infty\int_{F_\infty}^{F_n}{\bf 1}_A(x)dx.
\]
Therefore
\begin{eqnarray*}
&&\left|E\left[
\frac{\langle D\int_{F_\infty}^{F_n}{\bf 1}_A(x)dx,DF_\infty\rangle_\HH}{\|DF_\infty\|^2_\HH+\e}
\right]\right|\\
&=&
\left|E\left[
\int_{F_\infty}^{F_n}{\bf 1}_A(x)dx\left(
-\left\langle
DF_\infty,D\frac{1}{\|DF_\infty\|^2_\HH+\e}
\right\rangle_\HH
+\frac{LF_\infty}{\|DF_\infty\|^2_\HH+\e}
\right)
\right]\right|\\
&=&
\left|E\left[
\int_{F_\infty}^{F_n}{\bf 1}_A(x)dx\left(
\frac{2\langle D^2 F_\infty,D F_\infty\otimes D F_\infty\rangle_{\HH^{\otimes 2}}}{\left(\|DF_\infty\|^2_\HH+\e\right)^2}
+\frac{LF_\infty}{\|DF_\infty\|^2_\HH+\e}
\right)\right]\right|\leq\frac{c}{\e}\,\|F_n-F_\infty\|_{2},
\end{eqnarray*}
the last inequality following from Cauchy-Schwarz and the fact that $F_\infty\in\mathbb{D}^{2,4}$.
On the other hand,
\[
\left|E\left[{\bf 1}_A(F_n)\frac{\langle D(F_n-F_\infty),DF_\infty\rangle_\HH}{\|DF_\infty\|^2_\HH+\e}\right]\right|\leq
\frac{c}{\e}\|F_n-F_\infty\|_{\mathbb{D}^{1,2}}.
\]
Finally, let us observe that:
\begin{eqnarray*}
\left|E\left[\big({\bf 1}_A(F_n)-{\bf 1}_A(F_\infty)\big)\frac{\e}{\|DF_\infty\|_\HH^2+\e}\right]\right|
&\leq& E\left[\frac{\e}{\|DF_\infty\|_\HH^2+\e}\right]
\leq \e^\frac{\alpha}{2} \esp\left[\frac{1}{\|DF_\infty\|_\HH^\alpha}\right].
\end{eqnarray*}
Therefore, putting all these facts together and with $\e=\|F_n-F_\infty\|_{\mathbb{D}^{1,2}}^{\frac{2}{\alpha+2}}$, we get
\[
\big|P(F_n\in A)-P(F_\infty\in A)\big|\leq c\left(\e^\frac{\alpha}{2} +\frac{1}{\e}\|F_n-F_\infty\|_{\mathbb{D}^{1,2}}\right)\leq c \|F_n-F_\infty\|_{\D^{1,2}}^{\frac{\alpha}{\alpha+2}},
\]
which is the desired conclusion.
\qed

\bigskip
Let us now study the continuity of the law of $I_k(f)$ with respect to its kernel $f$. Before
offering another proof of the main result in Davydov and Martynova \cite{DM} (see our comments about this in the introduction), we start with a preliminary lemma.

\begin{lem}\label{prop-cw}
Let $F=I_k(f)$ with $k\geq 2$ and $f\in \HH^{\odot k}$ non identically zero.
There exists $c>0$ such that, for all $\e>0$,
\[
P(\|DF\|_\HH^2\leq \e)\leq c\,\e^{\frac1{2k-2}}.
\]
\end{lem}
{\it Proof}.
Throughout the proof, the letter $c$ stands for a non-negative constant {\it independent} of $n$
and whose value may change from line to line.
The proof is very close to that of Step 1 in Theorem \ref{Loi donne VT}.
Let $\{e_i\}_{i\geq 1}$ be an orthonormal basis of $\HH$. One can decompose
$f$ as
\begin{equation}\label{olala}
f=\sum_{i_1,\ldots,i_k=1}^\infty c_{i_1,\ldots,i_k}\,e_{i_1}\otimes\ldots\otimes e_{i_k}.
\end{equation}
For each $n\geq 1$, set
\[
f_n=\sum_{i_1,\ldots,i_k=1}^n c_{i_1,\ldots,i_k}\,e_{i_1}\otimes\ldots\otimes e_{i_k}.
\]
As $n\to\infty$, one has $f_n\to f$ in $\HH^{\otimes k}$ or, equivalently, $I_k(f_n)\to I_k(f)$ in $L^2(\Omega)$.
We deduce that there exists a strictly increasing sequence $\{n_l\}$ such that $I_k(f_{n_l})\to I_k(f)$ almost surely as $l\to\infty$.

On the other hand, this is a well-known result from It\^o that, with $k=k_1+\ldots+k_m$, one has
\[
I_k(e_1^{\otimes k_1}\otimes \ldots\otimes e_m^{\otimes k_m})=\prod_{i=1}^m H_{k_i}\left(X(e_i)\right),
\]
with $H_k$ the $k$th Hermite polynomial given by (\ref{hq}).
Also, one should note that the value of  $I_k(e_1^{\otimes k_1}\otimes \ldots\otimes e_m^{\otimes k_m})$
is not modified when one permutes the order of the elements in the tensor product.
It is deduced from these two facts that
\[
I_k(f_n)=Q_{n,k}\left(X(e_1),\ldots,X(e_n)\right),
\]
where $Q_{n,k}$ is a polynomial of degree at most $k$.
Theorem \ref{cw-thm}
ensures the existence of a constant $c>0$  such that,
for all $n\geq 1$ and $\e>0$,
\[
P\big(|I_k(f_n)|\leq \e\,\|f_n\|_{\HH^{\otimes k}}\big)\leq c\,\e^{1/k}.
\]
Next, we can use Fatou's lemma to deduce that, for any $\e>0$,
\begin{eqnarray*}
P\big(|I_k(f)|\leq \e\,\|f\|_{\HH^{\otimes k}}\big)&\leq&
P\left(\liminf_{l\to\infty} \{|I_k(f_{n_l})|\leq 2\e \|f_{n_l}\|_{\HH^{\otimes k}}\} \right)\\
&\leq&\liminf_{l\to\infty} P(|I_k(f_{n_l})|\leq 2\e \|f_{n_l}\|_{\HH^{\otimes k}})\leq c\,\e^{1/k}.
\end{eqnarray*}
Equivalently,
\[
P\big(|I_k(f)|\leq \e\big)\leq c\,\|f\|_{\HH^{\otimes k}}^{-1/k}\,\e^{1/k}.
\]
Now, assume for a while that $\langle f,h\rangle_{\HH} =0$
for all $h\in \HH$. By (\ref{olala}), we have
\[
\langle f,h\rangle_{\HH} =
\sum_{i_1,\ldots,i_k=1}^\infty c(i_1,\ldots,i_k)\langle e_{i_1},h\rangle_{\HH} \,\,\,e_{i_2}\otimes\ldots
\otimes e_{i_k},
\]
implying in turn, because $\langle f,h\rangle_{\HH}=0$ for all $h\in \HH$, that
\[
\sum_{i_2,\ldots,i_k=1}^\infty \left( \sum_{i_1=1}^\infty c(i_1,\ldots,i_k)\langle e_{i_1},h\rangle_{\HH}\right)^2=0\,\,\mbox{ for all $h\in \HH$.}
\]
By choosing $h=e_{i}$, $i=1,2,...$, we get that $c(i_1,\ldots,i_k)=0$ for any $i_1,\ldots,i_k\geq 1$, that is, $f=0$.
This latter fact being in contradiction with our assumption, one deduces that there exists $h\in \HH$
so that $\langle f,h\rangle_{\HH} \neq 0$.
Consequently,
\[
P\big(\|DF\|^2_\HH\leq \e\big)\leq
P\big(|\langle DF,h\rangle_\HH|\leq \sqrt{\e}\|h\|_\HH\big)=
P\left(|I_{k-1}(\langle f,h)\rangle_\HH|\leq \frac1k\sqrt{\e}\|h\|_\HH
\right)\leq c\,\e^{\frac1{2k-2}},
\]
which is the desired conclusion.
\qed

\bigskip

Finally, we state and prove the following result, which gives
a precise estimate for the continuity of $I_k(f)$ with respect to $f$.
This is almost the main result of Davydov and Martynova \cite{DM}, see our comments in the introduction.
Moreover, with respect to what we would have obtained by applying (\ref{Final}), here the rate is $\frac{1}{2k}$ (which
is better than $\frac{1}{2k+1}$, immediate consequence
of (\ref{Final})).

\begin{thm}\label{thm-dv}
Fix $k\geq 2$, and let $\{f_n\}_{n\geq 1}$ be a sequence of elements of $\HH^{\odot k}$.
Assume that $f_\infty=\lim_{n\to\infty}f_n$ exists in $\HH^{\otimes k}$ and that each $f_n$ as well as $f_\infty$ are not identically zero.
Then there exists a constant $c$, depending only on $k$ and $f_\infty$, such that, for all $n\geq 1$,
\[
d_{TV}(I_k(f_n),I_k(f_\infty))\leq c\,\|f_n-f_\infty\|_{\HH^{\otimes k}}^{\frac{1}{2k}}
\]
for any $n\geq 1$.
\end{thm}
{\it Proof}.
Set $F_n=I_k(f_n)$ and $F_\infty=I_k(f_\infty)$.
Let $A$ be a bounded Borel set of $\R$, and fix $0<\e\leq 1$.
Since $f_n,f_\infty\not\equiv 0$, Shigekawa theorem (see \cite{Shigekawa}, or \cite[Theorem 2.10.1]{NouPecBook}, or Theorem \ref{Loi donne VT}) ensures that $F_n$ and $F_\infty$ both have a density.
We deduce that
\begin{eqnarray*}
\left\langle D\int_{F_\infty}^{F_n} {\bf 1}_A(x)dx,DF_\infty\right\rangle_\HH &=&
\big({\bf 1}_A(F_n)-{\bf 1}_A(F_\infty)\big)\|DF_\infty\|_\HH^2+{\bf 1}_A(F_n)\langle D(F_n-F_\infty),DF_\infty\rangle_\HH,
\end{eqnarray*}
implying in turn that
\begin{eqnarray*}
P(F_n\in A)-P(F_\infty\in A)&=&
E\left[
\frac{\left\langle D\int_{F_\infty}^{F_n} {\bf 1}_A(x)dx,DF_\infty\right\rangle_\HH}{\|DF_\infty\|_\HH^2+\e}\right]\\
&&-E\left[\frac{{\bf 1}_A(F_n)\langle D(F_n-F_\infty),DF_\infty\rangle_\HH}{\|DF_\infty\|^2_\HH+\e}\right]\\
&&+E\left[\big({\bf 1}_A(F_n)-{\bf 1}_A(F_\infty)\big)\frac{\e}{\|DF_\infty\|_\HH^2+\e}\right].
\end{eqnarray*}
Firstly, using (\ref{etoile}) and next $\delta (DF_\infty)=-LF_\infty=kF_\infty$
we can write
\begin{eqnarray*}
&&\left|E\left[\frac{\left\langle D\int_{F_\infty}^{F_n} {\bf 1}_A(x)dx,DF_\infty\right\rangle_\HH}{\|DF_\infty\|^2_\HH+\e}\right]\right|\\
 &=&\left|E\left[
 \delta\left(\frac{DF_\infty}{\|DF_\infty\|^2_\HH+\e}\right)
 \int_{F_\infty}^{F_n} {\bf 1}_A(x)dx
 \right]\right|\\
  &=&\left|E\left[
 \frac{kF_\infty}{\|DF_\infty\|_\HH^2+\e}\, \int_{F_\infty}^{F_n} {\bf 1}_A(x)dx
 \right]
 -E\left[
 \left\langle DF_\infty,D\left(\frac{1}{\|DF_\infty\|_\HH^2+\e}\right)\right\rangle_\HH\, \int_{F_\infty}^{F_n} {\bf 1}_A(x)dx
 \right]\right|\\
   &=&\left|E\left[
 \frac{kF_\infty}{\|DF_\infty\|_\HH^2+\e}\, \int_{F_\infty}^{F_n} {\bf 1}_A(x)dx
 \right]
 +E\left[
 \frac{2\langle D^2F_\infty,DF_\infty\otimes DF_\infty\rangle_{\HH^{\otimes 2}}}{(\|DF_\infty\|_\HH^2+\e)^2}\, \int_{F_\infty}^{F_n} {\bf 1}_A(x)dx
 \right]\right|\\
    &\leq&\frac1\e\,E\left[
 \big(k|F_\infty|+2\|D^2F_\infty\|_{\HH^{\otimes 2}}\big)\left| F_n-F_\infty\right|
 \right]\\
 &\leq &\frac1\e\,\|f_n-f_\infty\|_{\HH^{\otimes k}}\,\sqrt{k!\,E\left[
 \big(k|F_\infty|+2\|D^2F_\infty\|_{\HH^{\otimes 2}}\big)^2\right]},
\end{eqnarray*}
where the last inequality comes from Cauchy-Schwarz and the isometry property of multiple integrals.
Secondly, using $\frac1kE\big[\|DF_\infty\|^2_\HH\big]=E\left[F_\infty\times \frac1k\delta DF_\infty\right]=E[F_\infty^2]$ and \[
\frac1kE\big[\|D(F_n-F_\infty)\|^2\big]=E[(F_n-F_\infty)^2]=k!\|f_n-f_\infty\|_{\HH^{\otimes k}}^2,
\]
we have
\begin{eqnarray*}
\left|E\left[\frac{{\bf 1}_A(F_n)\langle D(F_n-F_\infty),DF_\infty\rangle_\HH}{\|DF_\infty\|_\HH^2+\e}\right]\right|
\leq \frac1\e\,\|f_n-f_\infty\|_{\HH^{\otimes k}}\,\sqrt{k^2k!E[F_\infty^2]}.
\end{eqnarray*}
Thirdly,
\begin{eqnarray*}
&&\left|E\left[\big({\bf 1}_A(F_n)-{\bf 1}_A(F_\infty)\big)\frac{\e}{\|DF_\infty\|_\HH^2+\e}\right]\right|\leq
E\left[\frac{\e}{\|DF_\infty\|_\HH^2+\e}\right]\\
&\leq&E\left[\frac{\e}{\|DF_\infty\|_\HH^2+\e}\,\,{\bf 1}_{\left\{\|DF_\infty\|^2>\e^{\frac{2k-2}{2k-1}}\right\}}\right]
+P\left(\|DF_\infty\|_\HH^2\leq \e^{\frac{2k-2}{2k-1}}\right)\\
&\leq&\e^{\frac1{2k-1}}+c\,\e^{\frac1{2k-1}}=c\,\e^{\frac1{2k-1}},
\end{eqnarray*}
where the last inequality comes from Lemma \ref{prop-cw}.

By summarizing, we get
\[
\big|P(F_n\in A)-P(F_\infty\in A)\big|\leq \frac{c}{\e}\,\|f_n-f_\infty\|_{\HH^{\otimes k}}+ c\,\e^{\frac1{2k-1}}.
\]
The desired conclusion follows by choosing $\e=\|f_n-f_\infty\|_{\HH^{\otimes k}}^{\frac{2k-1}{2k}}$.
\qed

\section{The Peccati-Tudor theorem holds in total variation}\label{s:pt}

Let us first recall the Peccati-Tudor theorem \cite{PTu04}.

\begin{thm}\label{PecTud}
Let $d\geq 2$ and $k_d, \ldots, k_1\geq 1$ be some fixed
integers. Consider vectors
\[
F_n=(F_{1,n},\ldots,F_{d,n})=
(I_{k_1}(f_{1,n}),\ldots,I_{k_d}(f_{d,n})), \quad n\geq 1,
\]
with
$f_{i,n}\in \EuFrak{H}^{\odot k_i}$.
Let $N\sim\mathcal{N}_d(0,C)$ with $\det(C)>0$ and
assume that
\begin{equation}
\label{eq:asympcov}
\lim_{n\to\infty}
\esp[F_{i,n}F_{j,n}]=C(i,j),\quad 1\leq i,j\leq d.
\end{equation} Then, as $n\to\infty$,  the following two
conditions are equivalent:
\begin{itemize}
\item[(a)] $F_n$
converges in law to $N$;
\item[(b)] for every $1\leq i\leq d$, $F_{i,n}$
converges in law to $\mathcal{N}(0,C(i,i))$.
\end{itemize}
\end{thm}

The following result shows that the assertion $(a)$ in the previous theorem may be replaced for free by an a priori stronger assertion,
namely:
\begin{itemize}
\item[($a'$)] $F_n$
converges in total variation to $N$.
\end{itemize}
\begin{theo}\label{P.T}
Let $d\geq 2$ and $k_d, \ldots, k_1\geq 1$ be some fixed
integers. Consider vectors
\[
F_n=(F_{1,n},\ldots,F_{d,n})=
(I_{k_1}(f_{1,n}),\ldots,I_{k_d}(f_{d,n})), \quad n\geq 1,
\]
with
$f_{i,n}\in \EuFrak{H}^{\odot k_i}$.
As $n\to\infty$, assume that $F_n \overset{\rm law}{\to} N\sim\mathcal{N}_d(0,C)$ with $\det(C)>0$.
Then, $d_{TV}(F_n,N)\to 0$ as $n\to\infty$.
\end{theo}

During the proof of Theorem \ref{P.T}, we shall need the following auxiliary lemma.
(Recall that $\mathcal{H}_k$ denotes the $k$th Wiener chaos of $X$.)
\begin{lem}\label{lemme algèbre}
Let $\mathcal{A}$ be the class of sequences $\{Y_n\}_{n\geq 1}$ satisfying that:
$(i)$ there exists $p\in\N^*$ such that $Y_n\in\bigoplus_{k=0}^p\mathcal{H}_k$ for all $n$;
and $(ii)$ $\sup_{n\geq 1}\esp[Y_n^2]<\infty$.
 We have the following stability property for $\mathcal{A}$: if $\{Y_n\}_{n\geq 1}$ and $\{Z_n\}_{n\geq 1}$
 both belong to $\mathcal{A}$, then $\{\langle D Y_n,DZ_n\rangle_{\HH}\}_{n\geq 1}$ belongs to $\mathcal{A}$ too.
\end{lem}
{\it Proof}.
Let $\{Y_n\}_{n\geq 1}$ and $\{Z_n\}_{n\geq 1}$
 be two sequences of $\mathcal{A}$. We then have: $(i)$
 $Y_n=\esp[Y_n]+\sum_{k=1}^p I_k(g_{k,n})$ and $Z_n=\esp[Z_n]+\sum_{k=1}^p I_k(h_{k,n})$ for some integer $p$ and
 some elements $g_{k,n}$ and $h_{k,n}$ of $\HH^{\odot k}$; $(ii)$ $\sup_{n\geq 1}\|g_{k,n}\|_{\HH^{\otimes k}}^2<\infty$
 and $\sup_{n\geq 1}\|h_{k,n}\|_{\HH^{\otimes_k}}^2<\infty$ for all $k=1,\ldots,p$.
Using the product formula for multiple Wiener-It\^o integrals, it is straightforward to check that
\[
\langle DY_n,DZ_n\rangle_\HH=\sum_{k,l=1}^{p}kl\sum_{r=1}^{k\wedge l}(r-1)!\binom{k-1}{r-1}\binom{l-1}{r-1}I_{k+l-2r}(g_{k,n}\widetilde{\otimes}_r h_{l,n}).
\]
We deduce in particular that $\langle DY_n,DZ_n\rangle_\HH\in\bigoplus_{k=0}^{2p}\mathcal{H}_k$. Moreover,
since
\[
\|g_{k,n}\widetilde{\otimes}_r h_{l,n}\|_{\HH^{\otimes k+l-2r}}\leq
\|g_{k,n}\otimes_r h_{l,n}\|_{\HH^{\otimes k+l-2r}}\leq
\|g_{k,n}\|_{\HH^{\otimes k}}\| h_{l,n}\|_{\HH^{\otimes l}}\leq\frac12\big(\|g_{k,n}\|_{\HH^{\otimes k}}^2+\| h_{l,n}\|_{\HH^{\otimes l}}^2\big),
\]
we have that $\sup_{n\geq 1}\esp\big[\langle DY_n,DZ_n\rangle_\HH^2\big]<\infty$.
That is, the sequence $\{\langle DY_n,DZ_n\rangle_\HH\}_{n\geq 1}$ belongs to $\mathcal{A}$.
\qed

\bigskip

We are now in a position to prove Theorem \ref{P.T}.

\bigskip

\noindent{\it Proof of Theorem \ref{P.T}}.
First, using Lemma \ref{paley} and because $F_n\overset{\rm law}{\to}\mathcal{N}_d(0,C)$, it is straightforward to show
that $\esp[F_{i,n}F_{j,n}]\to C(i,j)$ as $n\to\infty$ for all $i,j=1,\ldots,d$.
Now, fix $M\geq 1$  and let $\phi\in\mathcal{C}^{\infty}_c([-M,M]^d)$. For any $i=1,\ldots,d$, define
\[
T_i[\phi](x)=\int_{0}^{x_i}\phi(x_1,\ldots,x_{i-1},t,x_{i+1},\ldots,x_d)dt,\quad x\in\R^d.
 \]
Also, set $T_{i_1,\ldots,i_a}=T_{i_1}\circ\ldots\circ T_{i_a}$, so that $\partial_{i_1,\ldots,i_a}T_{i_1,\ldots,i_a}[\phi]=\phi$.

\bigskip

The following lemma, which exhibits mere regularizing properties for the operators $T_i$, is going to play a crucial role in the proof.

\begin{lem}\label{lemme opérateur T}
The function $T_{d,\ldots,2,1}[\phi]$ satisfies the following two properties:
\begin{itemize}
\item  for all $k=1,\ldots,d$, $\|T_{k,\ldots,2,1}[\phi]\|_{\infty}\leq M^k\|\phi\|_{\infty}$;
\item for all $x,y\in\R^d$, $|T_{d,\ldots,2,1}[\phi](x)-T_{d,\ldots,2,1}[\phi](y)|\leq M^{d-1}\|\phi\|_{\infty}\|x-y\|_1$.
\end{itemize}
\end{lem}
{\it Proof}. For any $x,y\in\R^d$ and $k=1,\ldots,d$, we have
\begin{eqnarray*}
|T_{k,\ldots,2,1}[\phi](x)|
&\leq &\int_{0}^{|x_1|}\int_{0}^{|x_2|}\ldots\int_{0}^{|x_k|}|\phi(t_1,\ldots,t_k,x_{k+1},\ldots,x_d)|dt_1 dt_2 \ldots dt_k\\
&\leq &\int_{0}^{M}\int_{0}^{M}\ldots\int_{0}^{M}|\phi(t_1,\ldots,t_k,x_{k+1},\ldots,x_d)|dt_1 dt_2 \ldots d t_k\leq  M^k \|\phi\|_{\infty},
\end{eqnarray*}
whereas
\begin{eqnarray*}
&&|T_{1,2,\ldots,d}[\phi](x)-T_{1,2,\ldots,d}[\phi](y)|\\
&=&\left|\int_{0}^{x_1}\ldots \int_{0}^{x_d}\phi(t_1,\ldots,t_d)dt_1 dt_2 \ldots dt_d-\int_{0}^{y_1}\ldots\int_{0}^{y_d}\phi(t_1,\ldots,t_d)dt_1 dt_2 \ldots dt_d\right|\\
&=&\left|\sum_{i=1}^d\int_{0}^{x_1}\ldots \int_{0}^{x_{i-1}}\int_{x_i}^{y_{i}}\int_{0}^{y_{i+1}}\ldots\int_{0}^{y_d}\phi(t_1,\ldots,t_d)dt_1 dt_2 \ldots dt_d\right|
\leq M^{d-1}\|\phi\|_{\infty}\sum_{i=1}^d|x_i-y_i|.
\end{eqnarray*}
\qed

\bigskip

Let us go back to the proof of Theorem \ref{P.T}.
Set
\[
\Gamma_n= \left(
\begin{array}{ccc}
\langle D F_{1,n},D F_{1,n}\rangle_{\HH}&\ldots&\langle D F_{d,n},DF_{1,n}\rangle_{\HH}\\
\vdots&\ldots&\vdots\\
\langle D F_{1,n},D F_{d,n}\rangle_{\HH}&\ldots&\langle D F_{d,n},D F_{d,n}\rangle_{\HH}\\
\end{array}
\right),
\]
the Malliavin matrix associated with $F_n$.
Using the chain rule (\ref{e:chainrule}), we have
\begin{eqnarray}\label{système linéaire}
\left(
\begin{array}{c}
\langle D\phi(F_n),DF_{1,n}\rangle_{\HH}\\
\vdots\\
\langle D\phi(F_n),DF_{d,n}\rangle_{\HH}\\
\end{array}
\right)
&=&
\Gamma_n\left(
\begin{array}{c}
\partial_1\phi(F_n)\\
\vdots\\
\partial_d\phi(F_n)\\
\end{array}
\right).
\end{eqnarray}
Solving (\ref{système linéaire}) yields:
\begin{equation}\label{equation fondamentale}
\partial_i\phi(F_n)\det (\Gamma_n)=\sum_{a=1}^d ({\rm Com}\,\Gamma_n)_{a,i}\langle D\phi(F_n),DF_{a,n}\rangle_{\HH},
\end{equation}
where ${\rm Com}(\cdot)$ stands for the usual comatrice operator.
By first multiplying (\ref{equation fondamentale}) by $W\in\D^{1,2}$ and then taking the expectation, we get,
using (\ref{voiture}) as well,
\begin{eqnarray}\label{iter}
&&\esp[\partial_i\phi(F_n)\det (\Gamma_n)W]\\
&=&-\sum_{a=1}^d \esp\left[\phi(F_n)\big(\langle D(W({\rm Com} \,\Gamma_n)_{a,i}),DF_{a,n}\rangle_{\HH}+({\rm Com}\,\Gamma_n)_{a,i}W LF_{a,n}\big)\right]
=\esp[\phi(F_n)R_{i,n}(W)],\nonumber
\end{eqnarray}
where
\[
R_{i,n}(W)=-\sum_{a=1}^d\big(\langle D(W({\rm Com} \Gamma_n)_{a,i}),DF_{a,n}\rangle_{\HH}+({\rm Com}\Gamma_n)_{a,i}W LF_{a,n}\big).
\]
Thanks to \cite[Lemma 6]{NO}, we know that, for any $i,j=1,\ldots,d$,
\begin{equation}\label{p.t.Proof}
\langle D F_{i,n}, D F_{j,n}\rangle_{\HH}\overset{L^2}{\to} \sqrt{k_ik_j}\,C(i,j)\quad\mbox{as $n\to\infty$}.
\end{equation}
Also, Lemma \ref{lemme algèbre} implies that $\langle D F_{i,n}, DF_{j,n}\rangle_{\HH}$ is in a finite sum of chaoses and is bounded in $L^2(\Omega)$. By hypercontractivity, we deduce that $\langle D F_{i,n}, DF_{j,n}\rangle_{\HH}$ is actually bounded in all the $L^p(\Omega)$,
$p\geq 1$, and that the convergence in (\ref{p.t.Proof}) extends in all the $L^p(\Omega)$. As a consequence,
\begin{equation}\label{detcv}
\det(\Gamma_n)\overset{L^2}{\to} \det (C)\,\prod_{i=1}^d k_i=:\gamma>0.
\end{equation}
Using first (\ref{iter}) with $W=1$ and $T_1[\phi]$ instead of $\phi$, and then iterating, yields
\begin{eqnarray*}
\esp\left[\phi(F_n)\det(\Gamma_n)\right]&=&\esp\left[T_{1}[\phi](F_n)R_{1,n}(1)\right]\\
&=&\esp\left[T_{1}[\phi](F_n)\frac{\gamma-\det(\Gamma_n)}{\gamma}R_{1,n}(1)\right]+\esp\left[T_{2,1}[\phi](F_n)R_{2,n}\left(\frac{1}{\gamma} R_{1,n}(1)\right)\right]\\
&=&\ldots\\
&=&\sum_{k=1}^{d-1}\esp\left[T_{k,\ldots,1}[\phi](F_n)\frac{\gamma-\det(\Gamma_n)}{\gamma}P_{k,n}\right]+\esp\left[T_{d,\ldots,1}[\phi](F_n)P_{d,n}\right],
\end{eqnarray*}
with $P_{1,n}=R_{1,n}(1)$ and $P_{k+1,n}=R_{k+1,n}(\frac{1}{\gamma}P_{k,n})$.
As a consequence of Lemma \ref{lemme opérateur T}, we get the following inequality:
\begin{eqnarray}\label{fonda}
|\esp\left[\phi(F_n)\right]|&\leq&
\frac1\gamma \|\phi\|_\infty\,\|{\rm det}(\Gamma_n)-\gamma\|_{L^2}\\
&+&
\frac{\sum_{k=1}^{d-1} M^k\|\phi\|_{\infty}}{\gamma}\|\det(\Gamma_n)-\gamma\|_{L^2}\,\sup_{1\leq k\leq d-1}\|P_{k,n}\|_{L^2}+\|P_{d,n}\|_{L^2}\|T_{d,\ldots,1}[\phi]\|_{\infty}.\notag
\end{eqnarray}
Using Lemma \ref{lemme algèbre}, we have that $\{P_{k,n}\}_{n\geq 1}\in\mathcal{A}$ for all $k=1,\ldots,d$.
Hence, we arrive at the following inequality:
\[
|\esp\left[\phi(F_n)\right]|\leq c\left(\|\phi\|_{\infty}\|\det(\Gamma_n)-\gamma\|_{L^2}+\|T_{d,\ldots,1}[\phi]\|_{\infty}\right),\quad n\geq 1,
\]
where $c>0$ denote a constant independent of $n$, and whose value can freely change from line to line in what follows.
Similarly (more easily actually!), one also shows that
\[
|\esp\left[\phi(N)\right]|\leq c \|T_{d,\ldots,1}[\phi]\|_{\infty}.
\]
    Thus, if $\phi,\psi \in \mathcal{C}^{\infty}_c([-M,M]^d)$ are such that $\|\phi\|_{\infty}\leq 1$ and $\|\psi\|_{\infty}\leq 1$, we have,
    for all $n\geq 1$,
\begin{equation}\label{fondabisbisbis}
|\left(\esp\right[\phi(F_{n})\left]-\esp\left[\phi(N)\right]\right)-\left(\esp\left[\psi(F_{n})\right]-\esp\left[\psi(N)\right]\right)|\leq c \|\det\Gamma_{n}-\gamma\|_{L^2}+c\|T_{d,\ldots,1}[\phi-\psi]\|_{\infty}.
\end{equation}
Now, let $\rho : \R^d\longrightarrow \R_{+}$  be in $\mathcal{C}^{\infty}_c$ and satisfy $\int_{\R^d}\rho(x)dx=1$. As usual,
set $\rho_{\alpha}(x)=\frac{1}{\alpha^d}\rho(\frac{x}{\alpha})$ whenever $\alpha>0$.
\begin{lem}\label{utile}
For all $\alpha>0$, we have
\[
\|T_{d,\ldots,1}[\phi-\phi*\rho_{\alpha}]\|_{\infty}\leq 2 \alpha M^{d-1}\int_{\R^d}\|u\|_1\rho(u)du.
\]
\end{lem}
{\it Proof}.
We can write
\begin{eqnarray*}
&&T_{d,\ldots,1}[\phi-\phi*\rho_{\alpha}](x)\\
&=&\int_0^{x_1}\int_0^{x_2}\ldots\int_0^{x_d}\left[\int_{\R^d}\big(\phi(s_1,\ldots,s_d)-\phi(s_1-y_1,\ldots,s_d-y_d)\big)\rho_{\alpha}(y_1,\ldots,y_d)dy\right]ds\\
&=&\int_{\R^d}dy\,\rho_{\alpha}(y_1,\ldots,y_d)
\int_0^{x_1}\int_0^{x_2}\ldots\int_0^{x_d}\big(\phi(s_1,\ldots,s_d)-\phi(s_1-y_1,\ldots,s_d-y_d)\big)ds\\
&=&\int_{\R^d}dy\,\rho_{\alpha}(y_1,\ldots,y_d)\left(T_{d,\ldots,1}[\phi](x)-T_{d,\ldots,1}[\phi](x-y)-\int_{-y_1}^0\ldots\int_{-y_d}^0\phi(s)ds\right).
\end{eqnarray*}
According to Lemma \ref{lemme opérateur T}, we have
\[
|T_{d,\ldots,1}[\phi](x)-T_{d,\ldots,1}[\phi](x-y)|\leq M^{d-1}\|y\|_1\quad\mbox{and}\quad\left|\int_{-y_1}^0\ldots\int_{-y_d}^0\phi(s)ds\right|\leq M^{d-1}\|y\|_1.
\]
By combining these two bounds with the above equality, we get
\begin{eqnarray*}
|T_{d,\ldots,1}[\phi-\phi*\rho_{\alpha}](x)|\leq 2 M^{d-1}\int_{\R^d}\rho_{\alpha}(y)\|y\|_1dy
= 2\alpha M^{d-1}\int_{\R^d}\|u\|_1\rho(u)du,
\end{eqnarray*}
which is the announced result.
\qed

\bigskip

Using the previous lemma and applying (\ref{fondabisbisbis}) with $\psi=\phi*\rho_{\alpha}$, we deduce that, for some constant $c$
independent of $\alpha>0$ and $n\geq 1$,
\begin{equation}\label{gigatop}
|\left(\esp\left[\phi(F_{n})\right]-\esp\left[\phi(N)\right]\right)-\left(\esp\left[\phi*\rho_{\alpha}(F_{n})\right]-\esp\left[\phi*\rho_{\alpha}(N)\right]\right)|\leq
c \left(\|\det(\Gamma_{n})-\gamma\|_{L^2}+\alpha\right).
\end{equation}
But
\begin{eqnarray*}
|\phi*\rho_{\alpha}(x)-\phi*\rho_{\alpha}(x')|&\leq&\frac{1}{\alpha}\int_{\R^d}|\phi(y)|\left|\rho(\frac{x-y}{\alpha})-\rho(\frac{x'-y}{\alpha})\right|dy\\
&\leq&\frac{\|\rho'\|_\infty}{\alpha^2}\int_{\R^d}|\phi(y)|dy\,\,\|x-x'\|_1
\leq\frac{\|\rho'\|_\infty(2M)^{d}}{\alpha^2}\|x-x'\|_1,
\end{eqnarray*}
that is, $\phi*\rho_\alpha$ is Lipschitz continuous with a constant of the form $c/\alpha^2$.
We deduce that
\begin{equation}\label{inegz}
\left|\esp\left[\phi*\rho_{\alpha}(F_{n})\right]-\esp\left[\phi*\rho_{\alpha}(N)\right]\right|\leq\frac{c}{\alpha^2}d_{W}(F_{n},N),
\end{equation}
where $d_W(F_n,N)$ stands for the Wasserstein distance between $F_n$ and $N$, that is,
\[
d_W(F_n,N)=\sup_{\phi\in Lip(1)} \big| E[\phi(F_n)]-E[\phi(N)]\big|.
\]
By plugging inequality (\ref{inegz}) into (\ref{gigatop}), we deduce that, for all $\alpha>0$ and all $n\geq 1$,
\[
\sup_{\phi}|\esp\left[\phi(F_{n})\right]-\esp\left[\phi(N)\right]|\leq \frac{c}{\alpha^2}d_{W}(F_{n},N)+c \big(\|\det(\Gamma_n)-\gamma\|_{L^2} +  \alpha\big),
\]
where the supremum runs over the functions $\phi\in\mathcal{C}^{\infty}_c([-M,M]^d)$ with $\|\phi\|_{\infty}\leq 1$
and where $c$ is a constant independent of $n$ and $\alpha>0$.
By letting $n\to\infty$ (recall that $d_W(F_n,N)\to 0$ by \cite[Proposition 3.10]{multivariate} and
$\det(\Gamma_n)\overset{L^2}{\to}\gamma$ by (\ref{detcv})) and then $\alpha\to 0$, we get:
\[
\lim_{n\to\infty}\sup_{\phi}|\esp\left[\phi(F_{n})\right]-\esp\left[\phi(N)\right]|=0,
\]
so that the forthcoming Lemma \ref{tvreg} applies and allows to conclude.
\qed

\begin{lem}\label{tvreg}
Let $F_\infty$ and $F_n$ be random vectors of $\R^d$, $d\geq 1$. As $n\to\infty$, assume that
$F_n\overset{\rm law}{\to}F_\infty$ and that, for all $M\geq 1$,
\[
A_M(n):=\sup_{\phi}|\esp\left[\phi(F_{n})\right]-\esp\left[\phi(F_\infty)\right]|\to 0,
\]
where the supremum is taken over functions $\phi\in\mathcal{C}^{\infty}_c([-M,M]^d)$ which are bounded by 1.
Then $d_{TV}(F_n,F_\infty)\to 0$ as $n\to\infty$.
\end{lem}
{\it Proof}.
Let $\varepsilon>0$. Using the tightness of $F_n$, we get that there exists $M_{\varepsilon}$ large enough such that
\[\sup_n\bP(
\max_{1\leq i\leq d}|F_{i,n}|
\geq M_{\varepsilon})\leq\e\quad\mbox{and}\quad
\bP(
\max_{1\leq i\leq d}|F_{i,\infty}|
\geq M_{\varepsilon})\leq\varepsilon.\]
Let $\phi\in\mathcal{C}(\R^d,\R)$ with $\|\phi\|_{\infty}\leq 1$ and $M\geq M_{\varepsilon}+1$.
We have
\begin{eqnarray*}
\big|E\left[\phi(F_n)-\phi(F_{\infty})\right]\big|&\leq&\left|E\left[\mathbf{1}_{[-M_{\varepsilon},M_{\varepsilon}]^d}(F_n)
\phi(F_n)-\mathbf{1}_{[-M_{\varepsilon},M_{\varepsilon}]^d}(F_{\infty})\phi(F_{\infty})\right]\right|+2\varepsilon\\
&\leq&\sup_{\psi\in E_M}\big|E\left[\psi(F_n)-\psi(F_{\infty})\right]\big|+2\varepsilon
\leq A_M(n)+2\varepsilon.
\end{eqnarray*}
Here, $E_M$ is the set of smooth functions $\psi$ with compact support in $[-M,M]^d$ which are bounded by 1.
Hence, for all $\varepsilon>0$,
\[
\limsup_{n\to\infty}d_{TV}(F_n,F_{\infty})
=\frac12\,\limsup_{n\to\infty}\,\sup_{
\stackrel{
\phi\in\mathcal{C}(\R^d,\R):
}{
\|\phi\|_{\infty}\leq 1
}
} \big|E[\phi(F_n)- \phi(F_\infty)]\big|
\leq\varepsilon
\]
and the desired conclusion follows.
\qed
\bigskip

\noindent
{\bf Acknowledgements}. We would like to thank an anonymous referee for his/her very careful reading of the manuscript. Also, we are grateful to David Nualart for bringing his joint paper \cite{HuLuNu} to our attention.


\begin{thebibliography}{99}

\bibitem{BH}
\rm N. Bouleau and F. Hirsch (1991).
\it Dirichlet forms and analysis on Wiener space.
\rm W. de Gruyter Berlin.

\bibitem{breton}
\rm J.-C. Breton (2006).
\rm Convergence in variation of the joint laws of multiple Wiener-Itô integrals.
{\it Statist. Probab. Lett.} {\bf 76}, no. 17, 1904-1913.

\bibitem{CW}
\rm A. Carbery and J. Wright (2001):
\rm Distributional and $L^q$ norm inequalities for polynomials over convex bodies in $\R^n$.
\it Math. Research Lett. \rm {\bf 8}, 233-248.


\bibitem{ChatterjeeMeckes}
\rm S. Chatterjee and E. Meckes (2008):
\rm Multivariate normal approximation using exchangeable pairs.
\it ALEA \rm {\bf 4}, 257-283.

\bibitem{ChenGoldsteinShao} L. H. Y. Chen, L. Goldstein and Q.-M. Shao (2011). {\it Normal Approximation by Stein's Method}. Springer-Verlag, Berlin


\bibitem{DM} Y. A. Davydov and G. V. Martynova (1987).
Limit behavior of multiple stochastic integral.
{\it Statistics and control of random process}.
Preila, Nauka, Moscow, 55-57 (in Russian).



\bibitem{Dudley book} R.M.\ Dudley (2003). \textit{Real Analysis and
Probability }(2$^{\text{nd}}$ Edition). Cambridge University
Press, Cambridge.


\bibitem{HuLuNu}
Y. Hu, F. Lu and D. Nualart (2012).
\rm Convergence of densities of some nonlinear functionals of Gaussian processes.
Preprint.

\bibitem{L} M. Ledoux (2010).
\rm Chaos of a Markov operator and the fourth moment condition.
{\it Ann. Probab.}, to appear.


\bibitem{MalliavinHormander} P. Malliavin (1978). Stochastic calculus of variations and hypoelliptic operators. In: {\it Proc. Inter. Symp. on Stoch. Diff. Equations, Kyoto 1976},  195-263.

\bibitem{MOO}
E. Mossel, R. O'Donnell and K. Oleszkiewicz (2010). Noise stability of functions with low influences: invariance and optimality. {\it Ann. Math.} {\bf 171}, no. 1, 295-341.

\bibitem{E.N} E. Nelson (1973).
\rm The free Markoff field.
{\it J. Funct. Analysis} {\bf 12}, 211-227

\bibitem{webpage}
I. Nourdin.
\rm A special webpage on Stein's method and Malliavin calculus.
{\tt http://www.iecn.u-nancy.fr/}$\sim${\tt nourdin/steinmalliavin.htm}


\bibitem{lecturenotes-coursfondation} I. Nourdin (2012).
\rm Lectures on Gaussian approximation using Malliavin calculus.
{\tt http://www.iecn.u-nancy.fr/~nourdin/lecturenotes-coursfondation.pdf}


\bibitem{N-P}
\rm I. Nourdin and G. Peccati  (2009).
\rm Stein's method on Wiener chaos.
{\it Probab. Theory Rel. Fields} {\bf 145}(1), 75-118.


\bibitem{NouPecBook}
I. Nourdin and G. Peccati (2012).
{\it Normal Approximations Using Malliavin Calculus: from Stein's Method to Universality}.
Cambridge Tracts in Mathematics. Cambridge University Press.



\bibitem{multivariate}
\rm I. Nourdin, G. Peccati and A. R\'eveillac (2010).
\rm Multivariate normal approximation using Stein's method and Malliavin calculus.
{\it Ann. Inst. H. Poincar\'e (B) Probab. Statist.} {\bf 46}(1), 45-58.


\bibitem{G.I} I. Nourdin and G. Poly (2012).
\rm Convergence in law in the second Wiener/Wigner chaos.
{\it Electron. Comm. Probab.} {\bf 17}, no. 36.


\bibitem{nualartbook}
\rm D. Nualart (2006).
\it The Malliavin calculus and related topics of Probability and Its Applications.
\rm Springer-Verlag, Berlin, second edition.

\bibitem{NO}
D. Nualart and S. Ortiz-Latorre (2008).
\rm Central limit theorems for multiple stochastic integrals and Malliavin calculus.
\textit{Stoch. Proc. Appl.} {\bf 118} (4), 614-628.


\bibitem{nunugio}
\rm D. Nualart and G. Peccati (2005).
\rm Central limit theorems for sequences of multiple stochastic integrals.
{\it Ann. Probab.} {\bf 33} (1), 177-193.


\bibitem{PTu04} G. Peccati and C.A. Tudor (2005). Gaussian limits for
vector-valued multiple stochastic integrals. \textit{S\'{e}minaire de
Probabilit\'{e}s XXXVIII}, LNM \textbf{1857}. Springer-Verlag, pp. 247-262.

\bibitem{PM}
G. Poly and D. Malicet (2011).
\rm Properties of convergence in Dirichlet structures.
\rm Preprint.


\bibitem{Sch} M. Schreiber (1969). Fermeture en probabilit\'{e} de certains
sous-espaces d'un espace $L^{2}$. \textit{Zeitschrift Warsch. verw.\ Gebiete
}\textbf{14}, 36-48.

\bibitem{Shigekawa} I. Shigekawa (1980). Derivatives of Wiener functionals and absolute continuity of induced measures. {\it J. Math. Kyoto Univ.} {\bf 20}(2), 263-289.





\end{thebibliography}
\end{document}